\documentclass[11pt,tgrind,psbox,leqno]{article}
\usepackage{latexsym}
\usepackage{graphics}
\usepackage{amsmath}
\usepackage{xspace}
\usepackage{amssymb}

\newcommand {\be}[1]{\begin{equation}\label{#1}}
\newcommand {\ee}{\end{equation}}
\newcommand {\bea}{\begin{eqnarray}}
\newcommand {\eea}{\end{eqnarray}}
\newcommand{\Pois}{\ensuremath{\operatorname{Pois}}\xspace}

\newcommand{\qed}{\hfill $\Box$}
\newcommand{\remark}{{\bf Remark : }}

\newtheorem{theorem}{Theorem}
\newtheorem{lemma}[theorem]{Lemma}
\newtheorem{conj}{Conjecture}
\newtheorem{prop}{Proposition}

\newtheorem{Defi}{Definition}

\title{\bf Linear Phase Transition in  Random Linear Constraint Satisfaction Problems}

\author{David Gamarnik
\thanks{IBM T.J. Watson Research Center, Yorktown Heights, NY 10598, USA.
        Email address: gamarnik@watson.ibm.com }
}

\begin{document}

\maketitle

\centerline{{\bf Keywords:} Random K-SAT, Satisfiability Threshold, Linear Programming, Random Graphs}

\begin{abstract}
Our model is a generalized linear programming relaxation of a much studied random K-SAT problem.
Specifically, a set of linear constraints ${\cal C}$ on $K$ variables is fixed. From a pool of
$n$ variables, $K$ variables are chosen uniformly at random and a
constraint is chosen from ${\cal C}$ also uniformly at random.
This procedure is repeated $m$ times independently. We are interested
in whether the resulting linear programming problem
is feasible. We prove that the feasibility property experiences a
linear phase transition,  when $n\rightarrow\infty$ and $m=cn$ for
a constant $c$. Namely, there exists a critical value $c^*$
such that, when $c<c^*$, the problem is feasible or is
asymptotically almost feasible, as $n\rightarrow\infty$, but, when
$c>c^*$, the "distance" to feasibility is at least a positive
constant independent of $n$. Our result is obtained using the combination of a
powerful local weak convergence method developed in Aldous
\cite{Aldous:assignment92}, \cite{Aldous:assignment00}, Aldous
and Steele \cite{AldousSteele:survey}, Steele \cite{Steele:MinSpanningTree}
and martingale techniques.

By exploiting a linear programming duality, our theorem implies  the following result
in the context of sparse random graphs $G(n, cn)$ on $n$ nodes with
$cn$ edges, where edges are equipped with randomly generated weights.
Let ${\cal M}(n,c)$ denote maximum weight matching in $G(n, cn)$. We prove that when $c$ is a constant
and $n\rightarrow\infty$, the limit $\lim_{n\rightarrow\infty}{\cal M}(n,c)/n,$
exists, with high probability. We further extend this result to maximum weight $b$-matchings also
in $G(n,cn)$.
\end{abstract}

\section{Introduction}\label{introduction}

The primary objective  of the present paper is studying randomly generated linear programming
problems. We are interested in scaling behavior of the corresponding objective value and some phase transition properties,
as the size of the problem diverges to infinity.
Our random linear programming problems
are generated in a specific way. In particular, our linear programs
have a fixed number of variables per constraint and the number of variables and constraints
diverges to infinity in such a way that their ratio stays a constant.

Our motivation to consider this specific class of random linear programs has several
sources. The main motivation is recent explosion of interest
in random instances of boolean satisfiability (K-SAT) problems and ensuing phase transition phenomenon.
The main outstanding conjecture in this field states that the satisfiability property of random K-SAT problem experiences
a linear phase transition as the function of the ratio of the number of clauses to the number of variables. Our linear programming
problem can be viewed as a generalized linear programming relaxation of the integer programming
formulation of such random K-SAT problem.

Tightly related to the K-SAT  problem are problems
of maximal cardinality cuts, independent sets, matchings  and  other objects, in sparse random
graphs $G(n, \lfloor cn\rfloor)$, which are graphs on $n$ nodes and $ \lfloor cn\rfloor$ edges
selected uniformly at random from all the possible edges, and where  $c>0$ is  some fixed constant.
For future we drop the annoying notation $\lfloor\cdot\rfloor$, assuming that $cn$ is always an integer.
It is easy to show that all of these  objects scale linearly in $n$. It is conjectured that the size of each such object divided by
$n$ converges to a constant, independent of $n$. This convergence is established  only for the case
of maximal matchings using direct methods \cite{KarpSipser}, where the limit can be computed explicitly,
but is open in other cases.

The main result of this paper states that the objective value of the random linear
programming problem we consider, when divided by the number of variables converges with high probability (w.h.p.)
to a certain limit. As a corollary
we prove that, suitably defined, distance to feasibility in the same random linear programming problem
experiences a \emph{linear} phase transition, just as conjectured for random K-SAT problem.
Furthermore, we show that, in a special case, the dual of this random linear programming problem
is a linear programming relaxation of the  maximum cardinality matching and more generally $b$-matching
(defined later) problems in $G(n,  cn)$. We show that
these relaxations are asymptotically tight as the number of nodes $n$ diverges to
infinity. As a corollary of our main result, we prove that maximum cardinality $b$-matching
when divided by $n$ converge to a constant. These results hold even in the weighted version,
where edges are equipped with randomly independently generated non-negative weights.

Our  proof technique is a combination of a very powerful \emph{local weak convergence} method
and martingale techniques. The local weak convergence method was
developed in Aldous \cite{Aldous:assignment92}, \cite{Aldous:assignment00}, Aldous
and Steele \cite{AldousSteele:survey}, Steele \cite{Steele:MinSpanningTree}. The method
was specifically used by Aldous for proving the $\zeta(2)$ conjecture for the random
assignment problem. It was used in \cite{Aldous:assignment92} to prove
that the expected minimum weight matching in a complete bipartite graph  converges to a certain constant.
Later Aldous proved \cite{Aldous:assignment00}
that this limit is indeed $\zeta(2)$, as conjectured earlier by Mezard and Parisi \cite{MezardParisi}.
Since then the local weak convergence method was used for other problems (see \cite{AldousSteele:survey}
for a survey), and seems to be a very
useful method for proving existence of limits in  problems like the ones we described, and
in some instances also leads to actually computing the limits of interest. By an analogy with
the percolation literature, we call these problems existence and computation of \emph{scaling limits} in large random
combinatorial structures. Such questions, many of them open, abound in percolation literature.
For example the existence of limits of crossing probabilities in critical percolation have
been established in several percolation models like triangular percolation using conformal
invariance techniques \cite{SchrammPercolationFormula}, \cite{SmirnovWerner},
but are still open in the case of other lattices, like rectangular bond and site
critical percolation, see Langlands \cite{LanglandsPouliotSaint-Aubin}. Whether a local
weak convergence is a useful technique for addressing these questions seems worth investigation.

To the extend that we know, our result is the first application of the local weak convergence method to
establishing phase transitions in random combinatorial structures.
In the following section we describe in details randomly generated combinatorial problems we mentioned above,
describe the existing results in the literature and list some outstanding conjectures. In Section \ref{section:model}
we describe our model and state our main results. We also give a short summary of the proof steps. Sections
\ref{sec:PoissonTree}, \ref{sec:martingales}, \ref{section:projection} are devoted to the proof of
our main result. Section \ref{section:MatchingCyclePacking}, is devoted to the applications
to the maximum weight matching and $b$-matching in sparse random graphs. Section \ref{section:discussion} is
devoted to conclusions and some open problems.

\section{Background: random K-SAT, sparse random graphs and scaling limits}\label{subsection:ScalingLimits}

\subsection{Random K-SAT problem}\label{subsection:K-SAT}
A satisfiability or K-SAT problem is a boolean constraint satisfaction problem with a special form.
A collection of $n$ variables $x_1,x_2,\ldots,x_n$ with values in $\{0,1\}$ is fixed.
A boolean formulae of the form $C_1\wedge C_2\wedge \cdots \wedge C_m$ is constructed, where each $C_i$ is a disjunctive
clause of the form $x_{i_1}\vee \bar x_{i_2}\vee \bar x_{i_3}\vee\cdots \vee x_{i_K}$, where exactly $K$ variables are taken from the pool
$x_1,\ldots,x_n$, some with negation, some without.
The formulae is defined to be satisfiable if an assignment  of variables
$x_i, i=1,\ldots,n$ to $0$ or $1$ can be constructed such that all the clauses take value $1$.
The K-SAT problem is one of the most famous  combinatorial
optimization problem, see \cite{PapadimitriouBook}. It is well known that the satisfiability problem
is solvable in polynomial time for $K=2$, and is NP-complete for $K\geq 3$.

Recently we have witnessed an explosion of interest in random instances of the K-SAT problem.
This was motivated by computer science, artificial intelligence  and statistical physics investigations, with phase
transition phenomena becoming the focus of a particular attention. A random instance of a K-SAT problem with $m$ clauses and $n$
variables is  obtained by selecting each clause uniformly at random from the entire collection
of ${2^Kn^K\over K!}~$ ($2^K(\begin{array}{c}n \\ K\end{array})$) possible clauses where repetition of variables is (is not) allowed.
In particular, for each $j=1,2,\ldots,m,$
$K$ variables $x_{i_1},x_{i_2},\ldots,x_{i_K}$ or their negations are selected uniformly at random
from the pool $x_i,1\leq i\leq n$ to form a clause $C_j=y_{i_1}\vee\ldots y_{i_k}$, where
each $y_{i_r}=x_{i_r}$ or $\bar x_{i_r}$, equiprobably. This is done for all $j=1,\ldots,m$ independently. Whether the
resulting formulae has a satisfying assignment $\{x_1,\ldots,x_n\}\rightarrow \{0,1\}^n$
becomes a random event with respect to this random construction. The main outstanding conjecture
for the random K-SAT problem is as follows.

\begin{conj}\label{conj:main} For every $K\geq 2$ there exists a constant $c^*_K$ such that  a random K-SAT
formulae with $n$ variables and $m=cn$ clauses is
satisfiable when $c<c^*_K$ and  is not satisfiable when $c>c^*_K$, w.h.p. as $n\rightarrow\infty$. In other words,
the satisfiability experiences a \emph{linear} sharp  phase transition at $m=c^*_Kn$.
\end{conj}

That the problem experiences a sharp phase transition is proven by Friedghut \cite{Friedghut}
in  a much more general context. It is the linearity which is the main
outstanding feature of this conjecture. The conjecture can be rephrased as follows:
there does not exist $c_1>c_2$ and two infinite sequences $n^{(1)}_t,n^{(2)}_t,t=1,2,\ldots,$
such that instances of K-SAT problem with $n^{(1)}_t$ variables and $c_1n^{(1)}_t$ clauses
are satisfiable w.h.p., but instances  with $n^{(2)}_t$ variables and $c_2n^{(1)}_t$ clauses
are not satisfiable w.h.p., as $t\rightarrow\infty$. One of the goals
of our paper is to establish an analogue of this conjecture for generalized  linear programming relaxations of
the integer programming formulation (to be described below) of the random K-SAT problem.

Conjecture \ref{conj:main} is proven for the case $K=2$. Specifically,   $c^*_2=1$ was established by
Goerdt \cite{goerdt92}, \cite{goerdt96}, Chvatal and Reed \cite{ChvatalReed}, Fernandez de la Vega
\cite{vega}.
For higher values of $K$ many progressively sharper bounds on $c^*_K$ (assuming it exists) are established by now.
For $K=3$ the best known upper and lower bounds are $4.506$ and $3.42$,
obtained by Dubois, Boufkhad and Mandler \cite{DubBouMand}, and Kaporis, Kirousis and  Lalas \cite{KKL02}, respectively.
It is known that $c^*_K$, if exists, approaches asymptotically $2^K(\log 2+o(1))$ when $K$ is large, \cite{AchlioptasPeres}.
See also \cite{AchlioptasMoore}, \cite{FriezeKSATLargeK} for the related results. Talagrand \cite{TalagrandKSAT} approached
the random K-SAT problem using the methods of statistical physics.

The interest in random K-SAT problem does not stop at the threshold value $c^*_K$. For $c>c^*_K$
(assuming Conjecture \ref{conj:main} holds), a natural question is what is the maximal number of clauses $N(n,m)\leq m$ that
can be satisfied by a single assignment of the $n$ variables? It is shown in Coppersmith et al \cite{CopGamMohSor}
that for $K=2$ and  every $c>1$, there exists a constant $\alpha(c)>0$ such that $N(n,cn)\leq (c-\alpha(c))n$, w.h.p.
The following conjecture from (\cite{CopGamMohSor}) then   naturally extends  Conjecture \ref{conj:main}.

\begin{conj}\label{conj:MaxSAT}
Assuming Conjecture \ref{conj:main} holds,
$\lim_{n\rightarrow\infty}{N(n,cn)\over cn}$ exists and is smaller than one,
for all $c>c^*_K$.
\end{conj}

In Section \ref{section:model} we introduce a conjecture similar to the one above with respect to
random linear programs.

\subsection{Matching and $b$-matching  in $G(n, cn)$}\label{subsection:MatchingCyclePacking}
Let $G$ be a simple undirected graph on $n$ nodes $\{1,2,\ldots,n\}\equiv [n]$ with the edge set $E$.
A set of nodes $V\subset [n]$ in this graph is an independent set if no two nodes in $V$ are connected by an edge.
A partition of nodes $[n]$ into two $k$ groups $V_1,V_2,\ldots,V_k$ such that $\cup_{1\leq i\leq k} V_i=[n], V_{i_1}\cap V_{i_2}=\emptyset$
for all $i_1\neq i_2$, is defined to be a $k$-cut. The size of the  $k$-cut is the total number of edges whose end points
belong to different sets $V_i$. When $k=2$, the $k$-cuts are simply referred to as cuts.

A matching is a collection of edges such that no two edges are incident to the same node. The size of the
matching is the number of edges in it. A path is a collection of distinct nodes $C=\{i_1,\ldots,i_k\}$ such
that the edges $(i_1,i_2),(i_2,i_3),\ldots,$ $(i_{k-1},i_k)$ belong to the edge set $E$.
A cycle is a collection of distinct nodes $C=\{i_1,\ldots,i_k\}$ such that the edges $(i_1,i_2),(i_2,i_3),\ldots,$ $(i_{k-1},i_k),(i_k,i_1)$
belong to the edge set $E$.

Let $b\geq 1$ be a positive integer. A $b$-matching is a collection of edges $A\subset E$ such that
every node is incident to at most $b$ edges from $A$. Naturally, $1$-matching is simply a matching.
Note that $2$-matching is collection of node disjoint paths and cycles. We will also call it
path/cycle packing.

Fix a constant $c>0$. Let $G(n, cn)$ denote a simple undirected sparse random graph on $n$ nodes with $  cn$ edges selected
uniformly at random from all the possible $n(n-1)/2$ edges. This is a standard model of a sparse random graph.
Denote by ${\cal IND}(n,c),{\cal CUT}(n,c,k)$
${\cal M}(n,c,b)$ the size (carndinality) of the maximum independent set, $k$-cut and $b$-matching, respectively, in
$G(n,cn)$.  Suppose, in addition, the nodes and the edges of $G(n,cn)$ are
equipped with random non-negative weights $W^{\rm node}_{i,j},W^{\rm edge}_{i,j}$ drawn independently according to some common probability
distributions $\mathbb{P}\{W^{\rm node}\leq t\}\equiv w^{\rm node}(t), \mathbb{P}\{W^{\rm edge}\leq t\}\equiv w^{\rm edge}(t)$.
We assume throughout the paper that both $W^{\rm node}$ and $W^{\rm edge}$  have a bounded support $[0,B_w]$ (assumed the same for
simplicity). Let ${\cal IND}_w(n,c), {\cal CUT}_w(n,c,k)$
${\cal M}_w(n,c,b)$ denote the maximum  weight  independent set,  $k$-cut and $b$-matching,
respectively, where the weight of an independent set is the
sum of weights of its nodes, and the weights of a cut and $b$-matching are defined as
the sums of weights of edges in them.

It is well known and simple to prove that ${\cal IND}(n,c),$ ${\cal CUT}(n,c,k)$
${\cal M}(n,c,1)$ are all $\Theta(n)$  w.h.p. as $n$ diverges to infinity. For example
since a fixed node $i$ is isolated with a positive constant probability, then $\mathbb{E}[{\cal IND}(n,c)]=\Theta(n)$.
Since any matching is also a $b$-matching for $b\geq 1$, then ${\cal M}(n,c,b)=\Theta(n)$. Also the length of
the longest path in $G(n,cn)$ is also $\Theta(n)$, thanks to the result of Frieze \cite{FriezeLongPath}.

It is natural to suspect then that the expected values of these objects divided by $n$ converge to a constant,
both in the unweighted and weighted cases. In other words, the scaling limits exist for these objects.
In fact, it is conjectured in \cite{Aldous:FavoriteProblems} and
\cite{CopGamMohSor}, respectively, that the scaling limits
\be{eq:limits}
\lim_{n\rightarrow\infty}{\mathbb{E}[{\cal IND}(n,c)]\over n}, \qquad
\lim_{n\rightarrow\infty}{\mathbb{E}[{\cal CUT}(n,c)]\over n}
\ee
exist. The existence of these limits for expectation would also imply almost sure limits, by
application of Azuma's inequality.

The scaling limit of the form (\ref{eq:limits}) is in fact proven for maximum cardinality matchings ${\cal M}(n,c)$
 by Karp and Sipser \cite{KarpSipser}. The result was strengthened later by Aronson, Frieze and Pittel
\cite{AronsonFriezePittel}. Karp and Sipser proved that, almost surely,
\be{eq:KarpSipser}
\lim_{n\rightarrow\infty}{M(n,c)\over n}=1-{\gamma_*(c)+\gamma^*(c)+\gamma_*(c)\gamma^*(c)\over 2},
\ee
where $\gamma_*(c)$ is the smallest root of the equation $x=c\exp(-c\exp(-x))$ and $\gamma^*(c)=c\exp(-\gamma_*(c)$.
Their algorithmic method of proof is quite remarkable in its
simplicity. We briefly describe the argument below and explain why, however,  it does not
apply to the  case of weighted matchings. Suppose we are given a (non-random)
graph $G$. Then the following algorithm finds a maximum matching (clearly there could be many
maximum matchings): while the graph contains any leaves, pick any leaf of a tree and the corresponding edge. Mark the edge and delete the two
nodes incident  to the edge and all the other leaves that share the (deleted) parent with the selected leaf.
Delete all the edges incident to these leaves and delete the edge between the parent and its own parent. Repeat
while there are still leaves in the graph. When the graph has no more leaves select a maximum size matching
in the remaining graph. It is a fairly easy exercise to prove that this remaining  matching plus the set of marked edges is
a maximum matching. This fact is used to prove (\ref{eq:KarpSipser}).

One notes, however, that when edges of the graph  are equipped with some  weights, the Karp-Sipser algorithm
does not necessarily work anymore. Occasionally it might
be better to include an edge between a parent of a leaf and and a parent of a parent of a leaf and, as a result,
not include the edge incident to the leaf. Therefore, the Karp-Sipser algorithm may produce a strictly suboptimal
matching and the results (\ref{eq:KarpSipser}) do not hold for the weighted case. Moreover, it is not clear
how to extend the Karp-Sipser heuristic to $b$-matchings. In this paper we prove the convergence
(\ref{eq:KarpSipser}) for the maximum weight $b$-matchings. The proof uses the main result of the paper and
the linear programming duality, though  we are not able to compute the limits.
Naturally, our result applies to the non-weighted case -- maximum cardinality $b$-matching. To the best of
our knowledge this is a new result.

The case of maximum weight matching with random weights is treated by Aldous and Steele \cite{AldousSteele:survey} for the case of a
randomly generated tree on $n$ nodes. That is, consider a tree selected uniformly at random from the
set of all possible $n^{n-2}$ labelled trees. The limit of the sort (\ref{eq:KarpSipser}) is proven and computed
using the local weak convergence method, when the edges of this tree are equipped with exponentially
distributed random weights. The tree structure of the underlying graph helps very much the analysis.
In our case, however, the random graph $G(n, cn)$ contains a linear size non-tree "giant" component,
\cite{JansonBook}, when $c>1/2$, and the results of (\cite{AldousSteele:survey}) are not applicable.

Yet another scaling limit question is the existence of the limits for probability of $k$-colorability in $G(n, cn)$.
A graph is defined to be $k\geq 2$ colorable if there exists a function mapping vertices
of $G$ to colors $1,2,\ldots,k$ such that no two nodes connected by an edge have the same color.
The following conjecture proposed by Erdos is found in Alon and Spencer \cite{AlonSpencer}.

\begin{conj}\label{conj:coloring}
For every positive integer $k\geq 2$ there exists a critical value $c^*_k$ such
that the graph $G(n, cn)$ is w.h.p. $k$-colorable for $c<c^*_k$ and
w.h.p. not $k$-colorable for $c>c^*_k$.
\end{conj}
This conjecture is very similar in spirit to Conjecture \ref{conj:main}. For a survey
of existing results see a recent Molloy's survey \cite{MolloyColoringSurvey}. For related results
see also \cite{AchlioptasMooreColoring}, \cite{MooreMaxKCut}.

\section{Model and the main results}\label{section:model}
There is a natural way to describe a K-SAT problem as an integer programming problem.
The variables  are $x_i,i=1,2,\ldots,n$ which take values in $\{0,1\}$.
Each clause $C_j$ is replaced by a linear constraint of the form $x_{i_1}+(1-x_{i_2})+x_{i_3}+\ldots \geq 1$,
where term $(1-x)$ replaces $\bar x$.
For example a clause $C=x_3\vee x_7\vee \bar x_2\vee\bar x_4$ in a 4-SAT problem  is replaced by  a constraint
$x_3+x_7+(1-x_2)+(1-x_4)\geq 1$. It is easy to check that an assignment of $x_2,x_3,x_4,x_7$ to
$0$ and $1$ gives $C$ value $1$ if and only if the corresponding constraint is satisfied. Clearly,
these constraints can be created for all the possible clauses. In the present paper we study
the linear programming (LP) relaxation of this integer programming problem, where the restriction $x_i\in\{0,1\}$
is replaced by a weaker restriction $x_i\in [0,1]$. Note, that  this relaxation by itself is not
 interesting, as the assignment $x_i=1/2$ for all $i=1,2,\ldots,n$ makes all
of the linear constraints feasible. However, the problem becomes non-trivial when we generalize
the types of constraints that can be generated on the variables $x_i$, and this is described in the following subsection.

\subsection{Random K-LSAT problem}\label{subsection:mainresults}
Our setting is as follows. Consider a fixed collection of $K$ variables $y_1,y_2,\ldots,y_K$ which
take values in some bounded interval $B^1_x\leq y_i\leq B^2_x$ and a fixed collection ${\cal C}$ of
linear constraints on these variables: $\sum_{k=1}^K a_{rk}y_k\leq b_r, r=1,2,\ldots,|{\cal C}|$,
where the values $a_{rk},b_r$ are arbitrary  fixed reals. The $r$-th constraint can also be written
in  a vector form $a_ry\leq b_r$, where $a_r=(a_{r1},\ldots,a_{rK})$ and $y=(y_1,\ldots,y_K)$.
We fix $c>0$ and let $m=cn$, where $n$ is a large integer.
A random instance of a linear constraint satisfaction problem with $n+m$ variables
$x_1,\ldots,x_n, \psi_1,\ldots,\psi_m$
and $m$ constraints is constructed as follows. For each $j=1,2,\ldots,m$ we
perform the following operation independently. We
first select $K$ variables
$x_{i_1},x_{i_2},\ldots,x_{i_K}$ uniformly at random from $x_i,i=1,2,\ldots,n$.
Whether the variables are selected with or without replacement turns out to be irrelevant to the
results of this paper, as it is the case for random K-SAT problem. However, the order with which the variables are
selected  is relevant, since the constraints are not necessarily symmetric.
Then we select $1\leq r\leq |{\cal C}|$ also uniformly at random. We then generate a constraint

\be{eq:C}
C_j: \sum_{k=1}^Ka_{rk}x_{i_k}\leq b_r+\psi_j.
\ee
Here is an example of an instance with $K=3, n=10,m=4, |{\cal C}|=2$.
Say the first constraint $C_1$ is $2y_1+3y_2-y_3\leq 5,$ and the second constraint $C_2$ is
$-y_1+y_2+4y_3\leq 2$. An example of an instance where first three constraints are type $C_1$
and the fourth is type $C_2$ is
\begin{eqnarray*}
(2x_5+3x_4-x_9\leq 5+\psi_1)\wedge (2x_1+3x_3-x_4\leq 5+\psi_2) \wedge &  \\
(2x_2+3x_1-x_{10}\leq 5+\psi_3)  \wedge  (-x_5+x_8+4x_7\leq 2+\psi_4). &
\end{eqnarray*}

The central question is what are the optimal  values of $B^1_x\leq x_i\leq B^2_x, \psi_j\geq 0$,
which minimize the sum $\sum\psi_j$ subject to the constraints $C_j$. That is, we consider
the following linear programming problem:
\be{eq:LP}
{\rm Minimize}\sum_{1\leq j\leq m}\psi_j,\,\,
{\rm subject\,\,to:}\,\, C_1,C_2,\ldots,C_m,\,\,x_i\in [B^1_x,B^2_x], \psi_j\geq 0.
\ee

In words, we are seeking a solution $x_j$ which is as close to satisfying
the  constraints $\sum_{k=1}^Ka_{ri_k}x_{i_k}\leq b_r$ as possible. If the optimal value of this
linear programming problem is zero, that is $\psi_j=0$ for all $j$, then  all of these constraints
can be satisfied.
Naturally, the objective value of the linear program (\ref{eq:LP})   is a random variable.
We denote this random variable by ${\cal LP}(n,c)$. Note, that the linear program (\ref{eq:LP}) is
always feasible, by making $\psi_j$ sufficiently large. In fact, clearly, in the optimal solution
we must have $\psi_j=\max(0,\sum_{k=1}^Ka_{ri_k}x_{i_k}-b_r)$. We  refer to the  linear program (\ref{eq:LP})
as a \emph{random linear constraint satisfaction (LSAT) problem}, or \emph{random K-LSAT problem}.

The following conditions on the set of constraints ${\cal C}$ will be used  below.
\begin{itemize}
\item {\bf Condition ${\cal A}$.} For any constraint $a_ry\leq b_r, 1\leq r\leq |{\cal C}|$ for
any $k\leq K$ and for any value $z\in [B^1_x,B^2_x]$ there exist values $y_1,\ldots,y_K\in [B^1_x,B^2_x]$
such that $y_k=z$ and the constraint is satisfied.

\item{\bf Condition ${\cal B}$.} There exist a positive integer $l$ and a constant
$\nu>0$ such that for any $K$-dimensional cube $I$ of the form
$\prod_{1\leq k\leq K} [{i_k\over l},{i_k+1\over l}], B^1_x\leq {i_k\over l}< B^2_x$, $i_k$ integer, there
exists at least one constraint $\sum a_{rk}y_k\leq b_r$ from ${\cal C}$
such that for every $y\in I$, $\sum a_{rk}y_k-b_r\geq \nu$. That is, every point
of the cube $I$ deviates from satisfying this constraint by at least $\nu$.
\end{itemize}

The analogue of the Condition ${\cal A}$ clearly holds for random K-SAT problem. Given any clause
$y_1\vee y_2 \vee \cdots \vee y_K$ and $k\leq K$, if  $y_k$ is set to be $0$ or $1$, we still can
satisfy the clause, by satisfying any other variable. The following is an example of an LSAT problem where
Conditions ${\cal A}$ and ${\cal B}$ are satisfied.
Fix $K=3$. Let $B^1_x=0,B^2_x=1$, and let ${\cal C}$ be a collection of all eight constraints
of the type $-y_1-y_2-y_3\leq -7/4, -(1-y_1)-y_2-y_3\leq -7/4, \ldots, -(1-y_1)-(1-y_2)-(1-y_3)\leq -7/4$.
Condition ${\cal A}$ is checked trivially. We claim that Condition ${\cal B}$ holds for  $l=2$ and $\nu=1/4$.
Select any cube $I$ with side-length $1/l=1/2$. For example
$I=[0,1/2]\times [1/2,1]\times [1/2,1]$. Consider constraint $-y_1-(1-y_2)-(1-y_3)\leq -7/4$.
For any $y\in I$ we have $-y_1-(1-y_2)-(1-y_3)\geq -7/4+1/4=-7/4+\nu$. Other cases are
analyzed similarly.

Consider now the following generalization of the linear program (\ref{eq:LP}). For each $j=1,2,\ldots,m$
generate a random variable $W_j$, independently from some common distribution
$\mathbb{P}\{W_j\leq t\}$ with a bounded support $[-B_w,B_w]$. Let $w_x\geq 0$ and $w_{\psi}>0$ be  fixed non-negative
constants. Our random linear program in variables $x_i,\psi_j$ is constructed exactly as above except
each constraint $C_j:\sum_{1\leq r\leq K}a_{rk}x_{i_k}\leq b_r+\psi_j$ is replaced by
\be{eq:Cjgeneralized}
C_j:\sum_{1\leq r\leq K}a_{rk}x_{i_k}\leq b_r+W_j+\psi_j,
\ee
and the objective function is replaced by

\begin{eqnarray}
&{\rm Minimize}\,\,w_x\sum_{1\leq i\leq n}x_i+w_{\psi}\sum_{1\leq j\leq m}\psi_j,  \label{eq:GLP}\\
& \nonumber \\
&{\rm subject\,\,to:}\qquad C_1,C_2,\ldots,C_m,\,\,x_i\in [B^1_x,B^2_x],\,\, \psi_j\geq 0. \nonumber
\end{eqnarray}

This particular form of the linear program might look unnatural at first.
But note that setting $B_w=w_x=0,w_{\psi}=1$, turns this into exactly linear program (\ref{eq:LP}).
We will show later that this general format is useful when we study $b$-matchings in sparse
random graphs $G(n, cn)$. We denote the optimal value of the linear program (\ref{eq:GLP}) by
${\cal GLP}(n,c)$.  As before, this linear program is always feasible, by making $\psi_j$ sufficiently large.
Since we assumed $w_{\psi}>0$, then in the optimal solution
\be{eq:OptPsi}
\psi_j=\max(0,a_{rk}x_{i_k}-b_r-W_j).
\ee
We now state the main result of this paper. In words, our result asserts that the scaling limit of ${\cal GLP}(n,c)/n$ exists.

\begin{theorem}\label{theorem:mainMaxGen}
For every $c\geq 0$,  the limit
\be{eq:maxGLSAT}
\lim_{n\rightarrow\infty}{{\cal GLP}(n,c)\over n}\equiv f(c)
\ee
exists w.h.p. That is, there exists $f(c)\geq 0$ such that for every $\epsilon>0$,
\[
\mathbb{P}\{|{{\cal GLP}(n,c)\over n}-f(c)|>\epsilon\}\rightarrow 0
\]
as $n\rightarrow\infty$.
\end{theorem}

Our first application of Theorem \ref{theorem:mainMaxGen} is the following result.
It  establishes  a linear phase transition property
for the random K-LSAT problem. Recall that ${\cal LP}(n,c)$ is the
optimal value of the linear programming problem  (\ref{eq:LP}).
\begin{theorem}\label{theorem:mainMax}
There exists a constant $c^*_K>0$ such that, w.h.p. as $n\rightarrow\infty$,
\be{eq:MaxLSAT}
\lim_{n\rightarrow\infty}{{\cal LP}(n,c)\over n}=0,
\ee
for all $c<c^*_K$, and
\be{eq:MaxLSAT2}
\liminf_{n\rightarrow\infty}{{\cal LP}(n,c)\over n}>0,
\ee
for all $c>c^*_K$.
Moreover, if  Condition ${\cal A}$ holds, then $c^*_K>0$, and if Condition ${\cal B}$
holds, then $c^*_K<+\infty$.
\end{theorem}

In what sense does the theorem above establish a \emph{linear} phase transition?
It is conceivable that for a collection of constraints ${\cal C}$,
the following situation occurs:
there exist two constants $c_1>c_2$ and two sequences $n^{(1)}_t,n^{(2)}_t, t=0,1,2,\ldots$,
such that for $c=c_1$ the corresponding optimal values (\ref{eq:LP}) of the random K-LSAT problem
satisfy w.h.p. $\lim_t {\cal LP}(n^{(1)}_t,c)/n^{(1)}_t=0$, but for $c=c_2$,  $\liminf_t {\cal LP}(n^{(2)}_t,c)/n^{(2)}_t\geq \delta(c)>0$.
In other words, the critical density $c$ oscillates between different values.
This is precisely the behavior that Conjectures \ref{conj:main} and \ref{conj:MaxSAT} rule out for random K-SAT problem.
Our theorem states that such a thing is impossible for the random K-LSAT problem.
There exists a linear function
$c^*_Kn$ such that, w.h.p.,  below this function the instance is very close to being feasible, but
above this function the scaled "distance" $\min (1/n)\sum\psi_j$ to feasibility is at least
a positive constant.

The statement of the theorem above does not fully match its analogue, Conjecture \ref{conj:main},
as, using the auxiliary variables $\psi_j$ we converted the feasibility problem to the optimality problem.
Now consider the collection of constraints $C_j$ where $\psi_j$ are set to be zero, and we ask the
question whether the collection of constraints in (\ref{eq:LP}) has a feasible solution.
We suspect that this problem does experience a linear phase transition, but we do not have  a proof at the present time.

\begin{conj}\label{conj:LSAT}
Let $c^*_K$ be the value introduced in Theorem \ref{theorem:mainMax}. Then, w.h.p. as $n\rightarrow\infty$,
the random K-LSAT problem with $cn$ constraints is satisfiable if
$c<c^*_K$ and is not satisfiable if $c>c^*_K$.
\end{conj}

In this paper we use local weak convergence method  to prove Theorem \ref{theorem:mainMaxGen}.
While our  approach is very much similar to the one used in \cite{Aldous:assignment92}, there are several
distinctive features of our problem. In particular, we do not use an infinite tree construction
and instead consider a sequence of finite depth trees with some corresponding sequence of probability measures.
Then  we use a  Martingale Convergence Theorem for the  "projection" step.
This simplifies the proofs significantly.

\subsection{Maximum weighted $b$-matching}\label{subsection:mainresultsmatching}
We return to the setting of Subsection \ref{subsection:MatchingCyclePacking}. We have a sparse
random graph $G(n, cn)$, where $c$ is a positive constant. The edges of these graph are equipped
with random weights $W_{i,j}$ which are selected independently from a common distribution
$\mathbb{P}\{W_{i,j}\leq t\}=w^{\rm edge}(t), 0\leq t\leq B_w<\infty$, where $[0,B_w]$ is the support
of this distribution. Again let ${\cal M}_w(n,c,b)$
denote the maximum weight $b$-matching in $G(n, cn)$, where $b\geq 1$ is an integer.

\begin{theorem}\label{theorem:mainMatching}
For every $c>0$ the limit
\be{eq:LimitMatching}
\lim_{n\rightarrow\infty}{{\cal M}_w(n,c,b)\over n}\equiv g(c)
\ee
exists w.h.p.
\end{theorem}

The probability in the statement of the theorem is both with respect to the randomness of $G(n, cn)$ and
with respect to the random weights. This theorem is proven in Section \ref{section:MatchingCyclePacking}.
We use linear programming duality and certain linear programming formulation of the maximum weight $b$-matching problem
in order to related it to our main result, Theorem \ref{theorem:mainMaxGen}.

\subsection{Proof plan}\label{subsection:proofplan}
Below we outline the main steps in proving our main result, Theorem \ref{theorem:mainMaxGen}.
Let $\mathbb{E}[\cdot]$ denote the expectation operator.
The general scheme of the proof follows the one from Aldous \cite{Aldous:assignment92}.

\begin{enumerate}
\item We first observe that, as in the case of a random K-SAT problem, in the limit as $n\rightarrow\infty$,
the (random) number of constraints containing a fixed variable $x$ from the pool $x_1,\ldots,x_n$
is distributed as a Poisson random variable with parameter $cK$, denoted henceforth as $\Pois(cK)$.

\item For every $c>0$  we introduce
\be{eq:max}
\lambda(c)\equiv \liminf_{n\rightarrow\infty}{\mathbb{E}[{\cal GLP}(n,c)]\over n}.
\ee
Our goal is to show that in fact convergence $\lim_n{\mathbb{E}[{\cal GLP}(n,c)]\over n}$ holds,
and therefore we can set $f(c)=\lambda(c)$. The convergence w.h.p. will be a simple consequence
of Azuma's inequality. Then, in order to prove Theorem \ref{theorem:mainMax}, we prove that
$c^*_K\equiv \sup\{c:f(c)=0\}$ satisfies the properties required by the theorem.

\item We consider a subsequence $n_1,n_2,\ldots,n_i,\ldots$ along which
${\mathbb{E}[{\cal GLP}(n_i,c)]\over n_i}$ converges to $\lambda(c)$.
Let $X_1,\ldots,X_{n_i}$, $\Psi_1,\ldots,\Psi_{cn_i}\in [B^1_x,B^2_x]^{n_i}\times [0,\infty)^{cn_i}$
denote a (random) optimal assignment which
achieves the optimal value ${\cal GLP}(n_i)$.
For each $n_i$ we pick a variable $x_1$ from
the pool $x_1,\ldots,x_{n_i}$ (the actual index is irrelevant) and consider its depth $d$ neighborhood appropriately defined,
where $d$ is some fixed constant. We then consider the optimal solution $(X(n_i,d),\Psi(n_i,d))$
restricted to this $d$-neighborhood. We consider the  probability distribution
 ${\cal P}(d,n_i)$ which describes the joint probability distribution
for the values of $(X_i,\Psi_j,W_j)$ for $X_i,\Psi_j$ in the $d$-neighborhood
as well as the graph-theoretic structure of this neighborhood.

We show that for each fixed $d$, the sequence of probability measures ${\cal P}(d,n_i)$
is tight in its corresponding probability space.
As a result, there exists subsequence of $n_i$ along which the probability distribution ${\cal P}(d,n_i)$ converges
to a limiting probability distribution ${\cal P}(d)$ for every fixed $d$.
Moreover, we show that the subsequence can be selected in such a way that the
resulting probability distributions are consistent. Namely, for every $d<d'$, the marginal distribution
of ${\cal P}(d')$ in $d$-neighborhood is exactly ${\cal P}(d)$. We will show that, since the sequence $n_i$ was selected
to achieve the optimal value $\mathbb{E}[{\cal GLP}(n_i,c)]\approx \lambda(c)n$, then
\be{eq:locallimit}
\mathbb{E}[w_xX_1+{w_{\psi}\over K}\sum_j \Psi_j]=\lambda(c),
\ee
where the expectation is with respect to ${\cal P}(d)$ and the summation is over all the constraints $C_j$
containing $X_1$.

The sequence of probability distributions ${\cal P}(d), d=1,2,\ldots$ was used
by Aldous in \cite{Aldous:assignment92} to obtain an invariant (with respect to a certain pivot operator)
probability distribution on some infinite tree. Our proof does not require the analysis
of such a tree, although similar invariant measure can be constructed.

\item We consider a random sequence $\mathbb{E}[X_1|\Im_d],d=1,2,\ldots\,\,$, where $X_1\in [B^1_x,B^2_x]$
is, as above, the value that is assigned to the variable $x_1$ by an optimal solution, and
$\Im_d$ is the filtration corresponding to the sequence of probability measures ${\cal P}(d), d=1,2,\ldots\,\,$.
We prove that the sequence $\mathbb{E}[X_1|\Im_d],d=1,2,\ldots$ is a martingale.

\item This is the "projection" step in which for any $\epsilon>0$ and an arbitrary large $n$ we
construct a feasible solution to the system of constraints (\ref{eq:GLP})
which achieves the expected objective value at most $(\lambda(c)+\epsilon)n$. Given any large $n$ and an instance
of a random linear program (\ref{eq:GLP})
with variables $x_1,x_2,\ldots,x_n,\psi_1,\ldots,\psi_{cn}$
and  constraints $C_j, 1\leq j\leq cn$, for each variable $x_i,1\leq i\leq n$ we consider its $d$-neighborhood,
where $d$ is a very large constant. We let the value of $x_i$  be $\mathbb{E}[X_i|\Im_d]$
where the expectation is conditioned on the observed $d$-neighborhood of the variable $x_i$ and this
information is incorporated by filtration $\Im_d$. By construction,
this value is in $[B^1_x,B^2_x]$. Then we set $\Psi_j$ to the minimal value
which satisfies the $j$-th constraint for the selected values of $x_i$, for all $j=1,2,\ldots,cn$.
Using a martingale convergence theorem and property (\ref{eq:locallimit}) we show that
for a randomly chosen variable $x_i$, the corresponding value of  $\mathbb{E}[w_xX_i+{w_{\psi}\over K}\sum_j \Psi_j]$
is smaller than $\lambda(c)+\epsilon$, when $n$ and $d$ are sufficiently large. We sum the expectation above
over all $x_i$ and observe that each constraint belongs in the sum $\sum_j$ of exactly $K$ variables $x_i$.
Then the sum of these expectations is  $\mathbb{E}[w_x\sum_{1\leq i\leq n}X_i+w_{\psi}\sum_{1\leq j\leq cn}\Psi_j]$
which is exactly the objective function. We use this to conclude that the expected value of the objective
function is at most $(\lambda(c)+\epsilon)n$.
\end{enumerate}

\section{Poisson trees and some preliminary results}\label{sec:PoissonTree}
We begin by showing that in order to prove  Theorem \ref{theorem:mainMaxGen},
it suffices to prove the existence of a limit (\ref{eq:maxGLSAT}) for the
expected value of the optimal cost ${\cal GLP}(n,c)$. Indeed, note that
given an instance of a linear program (\ref{eq:GLP}), if we change one of the constraints $C_j$
to any other constraint from the pool ${\cal C}$ and change the value of $W_j$ to
any other value in $[-B_w,B_w]$, and leave all the other constraints intact, then the optimal value ${\cal GLP}(n,c)$
changes by at most a constant. Using a corollary of Azuma's inequality
(see Corollary 2.27 \cite{JansonBook} for the statement and a proof),
we obtain that $\mathbb{P}\{|{\cal GLP}(n,c)/n-\mathbb{E}[{\cal GLP}(n,c)/n]|>\epsilon\}$ converges to zero
exponentially fast for any $\epsilon>0$. Then the convergence $\lim_{n\rightarrow\infty}\mathbb{E}[{\cal GLP}(n,c)]/n$
implies the convergence of ${\cal GLP}(n,c)/n$ holds w.h.p.
Thus, from now on we concentrate on proving the existence of the limit
\be{eq:MaxLSATEXP}
\lim_n{\mathbb{E}[{\cal GLP}(n,cn)]\over n}.
\ee

A random instance  of a linear program (\ref{eq:GLP})
naturally leads to a sparse weighted $K$-hypergraph  structure on a node set $x_1,\ldots,x_n$.
Specifically, for each constraint $C_j,1\leq j\leq cn$ we create a $K$-edge
$(x_{i_1},\ldots,x_{i_K},r,w_j)$,
if $C_j$ contains exactly variables $x_{i_1},\ldots,x_{i_K}$ in this order,
the constraint is type $r, 1\leq r\leq |{\cal C}|$ and the corresponding random variable $W_j$ has value $w_j$.
This set of  edges   completely specifies the
random instance. We first study the distribution of the number of edges containing a
given fixed node $x=x_1,\ldots,x_n$.

\begin{lemma}\label{lemma:DegreePoisson}
Given node $x$ from the pool $x_1,\ldots,x_n$, the number of edges (constraints $C_j$) containing $x$
is distributed as $\Pois(cK)$, in the limit as $n\rightarrow\infty$.
\end{lemma}

\begin{proof} The probability that a given edge
does not contain $x$ is $1-K/n$ if variables are taken without replacement and
$((n-1)/n)^K=1-K/n+o(1/n)$ if taken with replacement. The probability
that exactly $s$ edges contain $x$ is then asymptotically $\Big(\begin{array}{c} cn \\s\end{array}\Big)
(K/n)^s(1-K/n)^{cn-s}$. When $s$
is a constant and $n\rightarrow\infty$, this converges to ${(cK)^s\over s!}\exp(-cK)$. \qed
\end{proof}

\vspace{.1in}

We now introduce a notion of a $d$-neighborhood of a variable $x$. A collection
of constraints $C_{i_1},C_{i_2},\ldots,C_{i_r}$, $1\leq i_j\leq m$ from an instance of  linear program (\ref{eq:GLP})
is defined to be a chain of length $r$ if for all $j=1,\ldots,r-1$ the constraints $C_{i_j}$ and $C_{i_{j+1}}$
share at least one variable. Fix a variable $x$ from the pool $x_i,1\leq i\leq n$.
We say that a variable $x'\in\{x_1,\ldots,x_n\}$ belongs to a $d$-neighborhood of $x$ if $x$ is connected to $x'$ by a chain
of length at most $d$. We say that a constraint $C_j$ belongs to the $d$-neighborhood of $x$ if
all of its  variables belong to it. The variables $W_j$ and $\psi_j$ in these constraints are
also assumed to be a part of this neighborhood.
In particular, a $1$-neighborhood of $x$ is the
set  of constraints $C_j$ which contain $x$ together with variables in these constraints.
If no constraints contain $x$, the $1$-neighborhood of $x$ is just $\{x\}$. For $d\geq 1$, we let
${\cal B}(x,d,n)$ denote the $d$-neighborhood of $x$, and let
$\partial {\cal B}(x,d,n)={\cal B}(x,d,n)\setminus {\cal B}(x,d-1,n)$ denote
the boundary of this neighborhood, where $\partial{\cal B}(x,0,n)$ is assumed to be $\emptyset$.
Of course ${\cal B}(x,d,n)$ and $\partial {\cal B}(x,d,n)$ are
random.
In graph-theoretic terms, ${\cal B}(x,d,n)$ is a sub-graph of the original $K$-hypergraph
corresponding to nodes with distance at most $d$ from $x$.


A chain $C_{i_1},C_{i_2},\ldots,C_{i_r}$ is defined to be a cycle if
$C_{i_1},C_{i_2},\ldots,C_{i_{r-1}}$ are distinct and $C_{i_r}=C_{i_1}$. The following
observation is a standard result from the theory of random graph. A simple proof is provided
for completeness.

\begin{lemma}\label{lemma:cycles}
Let $r,r'$ be  fixed constants. The expected number of cycles of length $r$ is at most $(K^2c)^r$.
Moreover, the expected number of variables with distance at most $r'$ from some size-$r$
cycle is at most $rc^{r+r'}K^{2r+r'}$.
In particular, for constants $r,r'$ the probability that a randomly and uniformly selected  variable  $x$ is at  distance at least $r'$ from any
size-$r$ cycle is at least $1-{rc^{r+r'}K^{2r+r'}\over n}$.
\end{lemma}

\begin{proof}
Fix any $r$ variables, say $x_1,\ldots,x_r$
and let us compute the expected number of cycles $C_1,\ldots,C_r$ such
that $C_j$ and $C_{j+1}$ share variable $x_j, j=1,\ldots,r$, (where $C_{r+1}$
is identified with $C_1$). Constraint $C_j$ contains variables $x_{j-1},x_j$
($x_0=x_r$). There are at most $n^{K-2}/(K-2)!$  choices for other variables
in $C_j$. For each such choice, the probability that a constraint consisting
of exactly this selection of variables is present in the random instance
is  at most  $(cn)/(n^K/K!)=K!c/n^{K-1}$. Finally, there are at most
$n^r$ ways to select (with order) $r$ variables $x_1,\ldots,x_r$. Combining, we obtain
that, asymptotically, the expected number of length-$r$ cycles is at most
\[
n^r\Big({n^{K-2}\over (K-2)!}K!{c\over n^{K-1}}\Big)^r< (K^2c)^r.
\]
By Lemma \ref{lemma:DegreePoisson}, for any fixed variable $x$ the expected number of variables with distance at most $r'$ from
$x$ is asymptotically at most $(cK)^{r'}$. Each size $r$ cycle contains at most $rK$ variables.
Therefore, the total expected number of variables with distance at most $r'$ to some
size $r$ cycle is at most $(K^2c)^r(rK)(cK)^{r'}=rc^{r+r'}K^{2r+r'}$. \qed
\end{proof}

\vspace{.1in}

Applying Lemmas \ref{lemma:DegreePoisson}, \ref{lemma:cycles} we obtain the following proposition,
which is well known in the context of random graphs \cite{JansonBook}.

\begin{prop}\label{prop:DegreePoisson}
Given a variable $x$, the number of constraints in ${\cal B}(x,1,n)$ is distributed as $\Pois(cK)$,
in the limit
as $n\rightarrow\infty$. Also ${\cal B}(x,d,n)$ is distributed as  a depth-$d$ Poisson tree. That is,
if $\partial{\cal B}(x,d,n)$
contains $r$ constraints and  $r(K-1)$ variables, then the number of constraints in
$\partial{\cal B}(x,d+1,n)$ is distributed as $\Pois(crK(K-1))$. Moreover, these constraints do not share variables,
other than variables in $\partial{\cal B}(x,d,n)$.
In short,  the constraints in ${\cal B}(x,d,n)$ are distributed as  first $d$ steps of a Galton-Watson (branching)
process with outdegree distribution $\Pois(cK)$, in the limit as $n\rightarrow\infty$.
\end{prop}

We finish this section by showing how Theorem \ref{theorem:mainMaxGen} implies Theorem \ref{theorem:mainMax}.
We noted before that the linear program (\ref{eq:LP}) is a special case of the linear program (\ref{eq:GLP}),
via setting $W_j=0$ with probability one and by setting $w_x=0, w_{\psi}=1$.
Assuming limit $f(c)$ defined in Theorem \ref{theorem:mainMaxGen} exists, let us first
show that it is a non-decreasing function of $c$. This is best seen by the following
coupling arguments. For any $c_1<c_2$ and large $n$ consider two instances of random
linear program (\ref{eq:GLP}) with $m=c_1n$ and $m=c_2n$, where the second is obtained by
adding $(c_2-c_1)n$ additional constraints to the first instance (we couple two instances). For each realization of
two linear programs, such an addition can only increase or leave the same the value of the
objective function. Note that in both cases we divide the objective value by the same $n$. We conclude
$f(c_1)\leq f(c_2)$.

Let $c^*_K\equiv \sup\{c\geq 0:f(c)=0\}$. Clearly the set in this definition includes $c=0$, and therefore is non-empty.
The definition of $c^*_K$ of course includes the possibilities  $c^*_K=0,c^*_K=\infty$ and clearly
satisfies (\ref{eq:MaxLSAT}) and (\ref{eq:MaxLSAT2}).

The proof of the second part of Theorem \ref{theorem:mainMax} follows from the following proposition.
\begin{prop}\label{prop:tildecK}
Under Condition ${\cal A}$, for any $c<1/K^2$,
a random K-LSAT instance has the optimal value
$\mathbb{E}[{\cal LP}(n,c)]=O(1)$.
Under Condition ${\cal B}$, there exists $\bar c_K>0$ such that  for all $c>\bar c_K$,
there exists $\delta(c)>0$ for which $\mathbb{E}[{\cal LP}(n,c)]\geq \delta(c)n$.
\end{prop}

\begin{proof}
Suppose Condition ${\cal A}$ holds. The technique we use is essentially
"set one variable at a time" algorithm, used in several papers on random K-SAT
problem to establish lower bounds on critical values of $c$, see, for example
\cite{AchlioptasSorkin}. Set variable $x_1$ to any value in $[B^1_x,B^2_x]$ (actual value is irrelevant).
For every constraint containing $x_1$ (if any exist), set other variables in this constraint so that this
constraints are satisfied with corresponding values of $\psi$ equal to zero.
This is possible by Condition ${\cal A}$. For every variable which is set
in this step, take all the other constraint containing this variable (if any exist)
and set its variables so that again these constraints are satisfied
with $\psi=0$. Continue in this fashion. If at any stage some newly considered
constraint $C_j$ contains some variable which was set in prior stages, set $\psi$
to the minimum value which guarantees to satisfy this constraint. At the worst we put $\psi=\max_r (|b_r|+\sum_k |a_{rk}|\max(|B^1_x|,|B^2_x|)$.
Note, however, that this situation occurs only if $C_j$ belongs to some cycle. Applying the first part of
Lemma \ref{lemma:cycles}, for $c<1/K^2$ the total expected number of constraints which belong to at least one cycle
is at most $\sum_{r=1}^{\infty} r(K^2c)^r=O(1)$. Therefore, the total expected number of constraints, for
which we set $\psi$ positive, is also $O(1)$ and $\mathbb{E}[{\cal LP}(n,c)]=O(1)$.

Suppose now Condition ${\cal B}$ holds. The proof is very similar to proofs of upper
bounds on critical thresholds for random K-SAT problems and uses the following moment argument.
Consider an instance of random K-LSAT with $n$ variables $m=cn$ constraints, where
$c$ is a very large constant, to be specified later.
Consider any $n$-dimensional cube
$I$ of the form $\prod_{1\leq k\leq n}[{i_k\over l},{i_{k+1}\over l}], B^1_x\leq {i_k\over l}<B^2_x$.
Consider the optimal cost ${\cal LP}(n,c)$ corresponding to solutions $x,\psi$, such
that $x$ is restricted to belong to $I$. Let, w.l.g.,  $x^j_1,\ldots,x^j_K$ be the
variables which belong to the $j$-th constraint, $1\leq j\leq m$. By Condition
${\cal B}$, with probability at least $1/|{\cal C}|$, $C_j$ is such that
the corresponding $\psi_j$ must be at least $\nu$. We obtain that the expected
cost $\mathbb{E}[{\cal LP}(n,c)]\geq \nu(1/|{\cal C}|)m$, when $x\in I$.
Moreover, since the events "$C_j$ is such a constraint" are independent for all $j$,
then, applying Chernoff bound, $\mathbb{P}\{{\cal LP}_I(n,c)\leq (1/2)\nu(1/|{\cal C}|)m\}\leq e^{-{m\over 8|{\cal C}|}}$
where the optimum ${\cal LP}_I(n,c)$ of the linear program (\ref{eq:LP}) is taken over the subset $x\in I$.
Since there are at most $l^n$ different cubes $I$, and every solution $x$ belong to at one of them
(several for points $x$ on the boundary between two cubes)
then $\mathbb{P}\{{\cal LP}(n,c)\leq (1/2)\nu(1/|{\cal C}|)m\}\leq l^ne^{-{m\over 8|{\cal C}|}}.$
By taking $m=cn$ with $c$ sufficiently large, we obtain that the probability above
is exponentially small in $n$. \qed
\end{proof}

\section{Limiting probability distributions}\label{section:limitingprobs}
In this and the following two sections we prove the existence of the limit (\ref{eq:MaxLSATEXP}).
Fix $c>0$ and take $n$ to be large.
We assume that the labelling (order) of the variables $x_1,\ldots,x_n$ is selected
independently from all the other randomness of the instance. In graph-theoretic
sense we have a random labelled hypergraph with labels of the nodes  independent from the
edges of the graph. We noted above that
\be{eq:maxPsi}
\psi_j\leq\max_{1\leq r\leq |{\cal C}|}\Big(\sum_k(|B^1_x|+|B^2_x|)K|a_{rk}|+|b_r|\Big)\equiv B_{\psi}.
\ee
Let, $X(n),\Psi(n)$ denote an optimal (random) assignment
of the random linear programming problem (\ref{eq:GLP}), where $X(n)=(X_1,\ldots,X_n),\Psi(n)=(\Psi_1,\ldots,\Psi_m)$.
That is $x_i=X_i,\psi_j=\Psi_j$ achieve the objective value ${\cal GLP}(n,c)$.
If the set of optimal solutions contains
more than one solution (in general it can be uncountable) select a solution $X(n),\Psi(n)$ uniformly
at random from this set. Define
\be{eq:lambdalimsup}
\lambda(c)=\liminf_n {\mathbb{E}[{\cal GLP}(n,c)]\over n},
\ee
and find a subsequence $n_t, t=1,2,\ldots$ along which
\be{eq:nt}
\lim_t{\mathbb{E}[{\cal GLP}(n_t,c)]\over n_t}=\lambda(c).
\ee
Fix a variable $x$ from the set of all $n_t$ $x$-variables.
Since labelling is chosen independent from the instance, we can
take $x=x_1$, w.l.g.
Denote by $X(d,n_t),\Psi(d,n_t),W(d,n_t)$ the collection of $X,\Psi$ and $W$-variables
which belong to the neighborhood ${\cal B}(x_1,d,n_t)$.
In particular $X_1\in X(d,n_t)$ and the number of $\Psi$ variables is the same as the number
of $W$ variables which is the number of constraints $C_j$ in ${\cal B}(x_1,d,n_t)$.
Denote by ${\cal P}(d,n_t)$ the joint probability distribution
of $({\cal B}(x_1,d,n_t),X(d,n_t),\Psi(d,n_t),W(d,n_t))$. We omit $x=x_1$ in the notation ${\cal P}(d,n_t)$ since, by symmetry,
this joint distribution is independent of the choice of $x$. The support of this probability distribution is
${\cal X}(d)\equiv\cup (T,[B^1_x,B^2_x]^{|T|}\times [0,B_{\psi}]^{E(T)}\times [-B_w,B_w]^{E(T)})\cup {\cal E}$,
where the first union runs over depth-$d$
rooted trees $T$ with root $x_1$, $|T|$ is number of nodes in $T$,
$E(T)$ is the number of constraints in $T$,
and ${\cal E}$ is a singleton event which represents  the event that ${\cal B}(x_1,d,n_t)$
is not a tree and contains a cycle. In particular, ${\cal X}(d)$ is a countable union of compact sets.
We have from Proposition \ref{prop:DegreePoisson}
that  $\lim_{t\rightarrow\infty}{\cal P}(d,n_t)({\cal E})=0$.
Observe that ${\cal X}(d)\subset {\cal X}(d+1)$ for all $d$.
As we showed in Proposition \ref{prop:DegreePoisson}, the marginal distribution of $T$ with respect to ${\cal P}(d,n_t)$
is depth-$d$ Poisson tree, in the limit as $t\rightarrow\infty$.

Recall, that a sequence of probability
measures ${\cal P}_n$ defined on a joint space ${\cal X}$ is said to be weakly converging to a probability
measure ${\cal P}$ if for any event $A\subset X$, $\lim_{n\rightarrow\infty}{\cal P}_n(A)={\cal P}(A)$.
We also need the following definition and a theorem, both can be found in \cite{durrett}.

\begin{Defi}\label{definition:tightness}
A sequence of probability measures ${\cal P}_n$ on ${\cal X}$ is defined
to be tight if for every $\epsilon>0$ there exists a compact set $K\subset {\cal X}$
such that $\limsup_n\mathbb{P}\{X\setminus K\}<\epsilon$.
\end{Defi}

\begin{theorem}\label{theorem:tightness}
Given a tight sequence of probability measures ${\cal P}_n$ on ${\cal X}$
there exists a probability measure ${\cal P}$ on ${\cal X}$ and a subsequence
${\cal P}_{n_t}$ that weakly converges to ${\cal P}$.
\end{theorem}

The following proposition is a key result of this section.

\begin{prop}\label{prop:ConvergenceMeasures}
For each $d=1,2,\ldots\,\,$ there exists a probability
measure ${\cal P}(d)$ on ${\cal X}(d)$
such that ${\cal P}(d,n_t)$ weakly converges to ${\cal P}(d)$. The sequence
of probability measures ${\cal P}(d), d=1,2,\ldots$ is consistent in the sense that
for every $d<d'$, ${\cal P}(d)$ is  the marginal distribution of ${\cal P}(d')$ on ${\cal X}(d)$.
The probability of the event ${\cal E}$ is equal to zero, with respect to ${\cal P}(d)$.
Finally, with respect to ${\cal P}(d)$,
\be{eq:lambdac}
\mathbb{E}[w_xX_1+{w_{\psi}\over K}\sum_j \Psi_j]=\lambda(c),
\ee
where the summation is over all the constraints $C_j$ in $T(d)$ containing the root variable $x_1$.
\end{prop}


\begin{proof} The proof is similar to the one in \cite{Aldous:assignment92} for constructing a
similar limiting probability distribution of optimal matchings in a complete bi-partite graphs
with random weights, where a compactness argument plus Kolmogorov's extension theorem is used to obtain the limiting
measures on infinite tree. In our case, since we limit ourselves to trees with bounded depths, the proofs can
be simplified by using the tightness argument. Fix $d\geq 1$. We claim that the
sequence of measures ${\cal P}(d,n_t)$ is tight on ${\cal X}(d)$.
 By Proposition \ref{prop:DegreePoisson}, according to the measure ${\cal P}(d,n_t)$ the neighborhood ${\cal B}(x_1,d,n_t)$
approaches in distribution a depth-$d$ Poisson tree with parameter $cK$. In particular,
the expected number of constraints in this neighborhood is smaller than $cK+c^2K^3+\ldots+c^dK^{2d-1}\equiv M_0$,
in the limit $t\rightarrow\infty$.
Fix $\epsilon>0$. By Markov's inequality the total number of constraints in ${\cal B}(x_1,d,n_t)$
is at most $M$ with probability at least $1-\epsilon$, for $M>M_0/\epsilon$ and $n_t$ sufficiently large.
This implies that, moreover, each $x$-variable in ${\cal B}(x_1,d,n_t)$ belongs to at most $M$ constraints
(has degree at most $M$) with probability at least $1-\epsilon$.

Let ${\cal X}(d,M)\subset {\cal X}(d)$ denote $\cup (T,[B^1_x,B^2_x]^{|T|}\times [0,B_{\psi}]^{E(T)}\times [-B_w,B_w]^{E(T)})$,
where the trees $T$ are restricted to have  degree
bounded by $M$.
That is,  each variable in such a tree belongs to at most $M$ constraints -- one
towards the root and $M-1$ in the opposite direction, and the root $x$ belongs to  at most $M$ constraints.
The number of trees  $T$ in ${\cal X}(d,M)$ is finite and, as a result, the set ${\cal X}(d,M)$
is compact, as it is a finite union of sets of the form $\{T\}\times [B^1_x,B^2_x]^{|T|}\times [0,B_{\psi}]^{|E(T)|}\times [-B_w,B_w]^{|E(T)|}$.
We showed above that according to ${\cal P}(d,n_t)$, the neighborhood ${\cal B}(x_1,d,n_t)$ belongs to ${\cal X}(d,M)$
with probability at least $1-\epsilon$, for all sufficiently large $t$. This proves that the sequence of measures ${\cal P}(d,n_t)$ is tight.
Then, applying Theorem \ref{theorem:tightness}, there exists a weakly converging subsequence ${\cal P}(d,n_{t_i})$.
We find such a sequence for $d=1$ and denote it by ${\cal P}(1,n^{(1)}_t), t=1,2,\ldots~$. Again using Theorem \ref{theorem:tightness},
for $d=2$ there exists a subsequence of ${\cal P}(2,n^{(1)}_t)$ which is weakly converging. We denote it by
${\cal P}(2,n^{(2)}_t)$. We continue this for all $d$ obtaining a chain of sequences
$n^{(1)}_t\supset n^{(2)}_t\supset\cdots \supset n^{(d)}_t\supset\cdots\,\,.$ Select a diagonal
subsequence $n^*_t\equiv n^{(t)}_t$ from these sequences. Then all the convergence above holds for this diagonal subsequence.
In particular, for every $d$, ${\cal P}(d,n^*_t)$ converges to some probability measure ${\cal P}(d)$.
Moreover, these measures, by construction, are consistent. Meaning, for every $d<d'$, the marginal
distribution of ${\cal P}(d')$ onto ${\cal X}(d)$ is simply ${\cal P}(d)$. Since, from
Proposition \ref{prop:DegreePoisson}, the probability of the event ${\cal E}$ according to ${\cal P}(d,n_t)$
approaches zero, then the probability of the same event with respect to ${\cal P}(d)$ is just zero, for every $d$.

To complete the proof of the proposition we need to establish (\ref{eq:lambdac}). Note that in
random instances of linear program (\ref{eq:GLP}), with $n_t$ $x$-variables, when we sum the expression
$w_xX_i+{w_{\psi}\over K}\sum_j\Psi_j$ over all  $x_i, i=1,2,\ldots,n_t$ we obtain
$w_x\sum_{1\leq i\leq n_t} X_i+w_{\psi}\sum_{1\leq j\leq cn_t}\Psi_j$, since each
variable  $\Psi_j$ appeared in exactly $K$ sums, corresponding to $K$ variables in $j$-th
constraint. From (\ref{eq:nt}) we obtain that $\lim_t\mathbb{E}[w_xX_i+{w_{\psi}\over K}\sum_j\Psi_j]=\lambda(c)$,
where the expectation is with respect to measure ${\cal P}(d,n_t)$. Since all the random variables involved
in $w_xX_i+{w_{\psi}\over K}\sum_j\Psi_j$ are bounded and since ${\cal P}(d,n^*_t)$ converges weakly to ${\cal P}(d)$,
then the convergence carries on to the expectations. This implies that with respect to ${\cal P}(d)$,
$\mathbb{E}[w_xX_i+{w_{\psi}\over K}\sum_j \Psi_j]=\lambda(c)$. \qed
\end{proof}

\section{Filtration and a martingale with respect to ${\cal P}(d)$}\label{sec:martingales}
Let us summarize the results of the previous section. We considered a sequence of
spaces ${\cal X}(1)\subset{\cal X}(2)\subset{\cal X}(3)\subset \cdots,$ where
${\cal X}(d)=\cup (T,[B^1_x,B^2_x]^{|T|}\times [0,B_{\psi}]^{|E(T)|}\times [-B_w,B_w]^{|E(T)|})$ and the union runs over
all depth-$d$ trees $T$. Recall, that the event ${\cal E}$ was used before to generically
represent the event that ${\cal B}(x,d,n)$ is not a tree. The probability of this event
is zero with respect to the limiting probability distribution ${\cal P}(d)$ we constructed,
so we drop this event from the space ${\cal X}(d)$.
We have constructed a probability measure ${\cal P}(d)$ on each ${\cal X}(d)$ as a weak limit of probability
measures ${\cal P}(d,n_t^*)$, defined on  $d$-neighborhoods of a  variable $x_1$. Denote by
$(T(d),X(d),\Psi(d),W(d))$ the random vector distributed according to ${\cal P}(d)$. Note, that $T(d)$
and $W(d)$ are independent from each other, first distributed as a depth-$d$ Poisson tree, second distributed
as i.i.d. with the distribution function $\mathbb{P}\{W\leq t\}$, yet $X(d)$ and $\Psi(d)$ depend on both $T(d)$ and $W(d)$.

Using Kolmogorov's extension theorem, the probability measures ${\cal P}(d)$ can be extended  to a probability
measure ${\cal P}$ on space $\cup (T,[B^1_x,B^2_x]^{\infty}\times [0,B_{\psi}]^{\infty})\times [-B_w,B_w]^{\infty})$, where the union runs
over all \emph{finite and infinite} trees $T$. This extension is used in
Aldous \cite{Aldous:assignment92} to construct a spatially invariant measure on infinite trees.
For our purpose this extension, while possible, is not necessary, and instead we use
a martingale convergence techniques.

Note, that, since the measures ${\cal P}(d)$ are consistent, they correspond to a certain filtration
on the increasing sequence of sets ${\cal X}(d),d=1,2,\ldots\,\,$.
Then we can look at $\mathbb{E}[X_1|(T(d),W(d))]$ as a stochastic process indexed by $d=1,2,\ldots\,\,$.
Another way of looking at this stochastic process is as follows.  A root $x=x_1$ is fixed. At time $d=1$ we
sample  constraints $C_j$ containing $x_1$ with the corresponding values of $W_j$ using  the ${\cal P}(1)$ distribution.
Recall that the number of constraints is distributed as $\Pois(cK)$, the type of each constraint is chosen
uniformly at random from $|{\cal C}|$ types, and $W_j$ values are selected i.i.d. using distribution function
$\mathbb{P}\{W\leq t\}$. For this sample we have a conditional probability distribution of $(X(1),\Psi(1))$,
that is of $x$ and $\psi$ variables in $1$-neighborhood of $x_1$.
Then, for each variable  in this depth-$1$
tree except for $x_1$ we again sample  constraints and $W_j$ values
to obtain a depth-$2$ Poisson tree $T(2)$ and $W(2)$. We obtain a distribution of $(X(2),\Psi(2))$ conditioned
on $T(2),W(2)$, and so on. On every step we reveal a deeper layer of the tree and and the corresponding values of $W_j$-s to
obtain a new conditional probability distribution for $(X(d),\Psi(d))$.

The  technical lemma that follows is needed to prove Proposition \ref{prop:martingale1} below.
This lemma is sometimes defined as \emph{tower property} of conditional expectations.
The proof of it can be found in many books on probability and measure theory, see for example
Theorem 1.2, Chapter 4 in \cite{durrett}.

\begin{lemma}\label{lemma:conditioning}
Let $X$ and $Y$ be dependent in general random variables, where $X\in R$ and $Y$ takes values in some general set
${\cal Y}$.  Let $f:{\cal Y}\rightarrow {\cal Y}$ be
a measurable function. Then $\mathbb{E}[\mathbb{E}[X|Y]|f(Y)]=\mathbb{E}[X|f(Y)]$.
\end{lemma}

The conditional expectations in this lemma are understood as follows. $\mathbb{E}[X|Y]$ is
an expectation of $X$ conditioned on the $\sigma_Y$-algebra generated by the random variable $Y$.
Similarly, $\mathbb{E}[X|f(Y)]$ is an expectation of $X$ conditioned on the $\sigma_{f(Y)}$-algebra generated
by the random variable $f(Y)$. Of course $\sigma_{f(Y)}\subset\sigma_Y$.

As above let $x_1$ be the root of our random tree $T(d)$ and let $X_1$
be the corresponding random variable from the vector  $(T(d),X(d),\Psi(d),W(d))$.
Conditioned on the event that $x_1$ belongs to at least one constraint $C_j$,
select any such a constraint and let, w.l.g. $x_2,\ldots,x_K$ be the variables in this constraint.
For every $k=2,\ldots,K$ and every $d\geq 2$, consider ${\cal B}(x_k,d)\equiv T(x_k,d)$, the
$d$-neighborhood of the variable $x_k$. One way to introduce this neighborhood formally is to select any $d'>d+1$, sample
$T(d')$ from ${\cal P}(d')$ and then consider the $d$-neighborhood around $x_k$, as a subtree
of the tree $T(d')$. But since, by consistency,  the distribution of these neighborhood is the same for all $d'>d+1$,
then we can simply speak about selecting $T(x_k,d)$. Note, that $T(x_1,d)=T(d)$. Let $W(x_k,d)$ denote the
collection of $W_j$ variables corresponding to constraints in $T(x_k,d)$. To simplify the notations, we
will let $T^W(x_k,d)$ stand for the pair $(T(x_k,d),W(x_k,d)$).

\begin{prop}\label{prop:martingale1}
For every $k=1,2,\ldots,K$ a random sequence $\mathbb{E}[X_k|T^W(x_k,d)],\, d=1,2,\ldots$ is a martingale
with values in $[B^1_x,B^2_x]$. Moreover, the sequence $\mathbb{E}[X_k|T^W(x_{d\mod(K)},d)],\, d=1,2,\ldots$
is also a martingale.
\end{prop}

\remark  To understand better the meaning of the first part of the proposition, imagine
that we first sample depth-$1$ tree $T^W(1)$. In case this tree is not trivial
($T^W(1)\neq \{x_1\}$), we fix a constraint from this tree and variable $x_k$ from this constraint.
Then we start revealing trees $T^W(x_k,2),T^W(x_k,3),\ldots\,.$ This defines the stochastic
process of interest.

The second part states that even if we reveal trees rooted in a round-robin fashion at variables $x_1,x_2,\ldots,x_K,x_1,x_2,\ldots,$
then we still obtain a martingale. That is the sequence $\mathbb{E}[X_k|T^W(x_1,1)]$, $\mathbb{E}[X_k|T^W(x_2,2)],\ldots,$
$\mathbb{E}[X_k|T^W(x_K,K)]$, $\mathbb{E}[X_k|T^W(x_1,K+1)]$, $\mathbb{E}[X_k|T^W(x_2,K+2)],\ldots,\,\,$ is a martingale.

\begin{proof} Recall, that the optimal values $X(d)$ are a weak limit
of optimal values $X(d,n_t^*), t=1,2,\ldots\,\,$. Then,
every  $X_i\in [B^1_x,B^2_x]$ almost surely.
To prove the martingale property we use Lemma \ref{lemma:conditioning}, where
$X=X_k, Y=T^W(x_k,d+1)$ and $f$  is a projection function which projects a depth-$d+1$
tree $T^W(x_k,d+1)$ onto a depth-$d$ tree $T^W(x_k,d)$ by truncating the $d+1$-st layer.  Applying the lemma, we obtain
$\mathbb{E}[\mathbb{E}[X_k|T^W(x_k,d+1)]|T^W(d)]=\mathbb{E}[X_k|T^W(x_k,d)]$, which means precisely that the
sequence is a martingale. The proof of the second part is exactly the same, we just observe that
$T^W(x_1,1)\subset T^W(x_2,2)\subset\cdots\subset T^W(x_{d\mod(K)},d)\subset\cdots\,\,,$
and we let $f$ to be a projection of $T^W(x_{d+1 \mod(K)},d+1)$ onto $T^W(x_{d\mod(K)},d)$. \qed
\end{proof}

The following is a classical result from the probability theory.

\begin{theorem}\label{theorem:MartingaleConvergence} {\bf [Martingale Convergence Theorem]}.
Let $X_n,n=1,2,\ldots\,\,$ be a martingale which takes values in some finite interval
$[-B,B]$. Then $X_n$ converges almost surely to some random variable $X$.
As a consequence, $|X_{n+1}-X_n|$ converges to zero almost surely.
\end{theorem}

A proof of this fact can be found in \cite{durrett}.
Here we provide, for completeness, a very simple proof of a weak version of MCT,
stating that $|X_{n+1}-X_n|$ converges to zero in probability. This is the only version we  need
for this paper. We have
\be{eq:martingale3}
0\leq \mathbb{E}[(X_{n+1}-X_n)^2]=\mathbb{E}[X^2_{n+1}]-2\mathbb{E}[X_{n+1}X_n]+\mathbb{E}[X^2_n].
\ee
But $\mathbb{E}[X_{n+1}X_n]=\mathbb{E}[(\mathbb{E}[X_{n+1}|X_n])X_n]=\mathbb{E}[X^2_n]$, where $\mathbb{E}[X_{n+1}|X_n]=X_n$ is
used in the second equality. Combining with (\ref{eq:martingale3}),
we obtain $0\leq \mathbb{E}[(X_{n+1}-X_n)^2]\leq \mathbb{E}[X^2_{n+1}]-\mathbb{E}[X^2_n]$ and therefore $\mathbb{E}[X_n^2]$ is an increasing subsequence.
Since $\mathbb{E}[X_n^2]\leq B^2$, then this sequence is converging and, from (\ref{eq:martingale3}),
$\mathbb{E}[(X_{n+1}-X_n)^2]$ converges to zero. Applying Markov's bound, we obtain that $|X_{n+1}-X_n|$ converges to zero in probability.

Let ${\cal J}(x_1)$ denote the set of constraints $C_j$ in $T(d)$ containing $x_1$ and let $J(x_1)=|{\cal J}(x_1)|$.
Select again any constraint from ${\cal J}(x_1)$. We denote this constraint by $C(x_1)$
and let it be $\sum_{1\leq k\leq K} a_{rk}y_k\leq b_r+W_j+\psi_j$ for some $1\leq r\leq |{\cal C}|$, in an expanded form.
Again let $x_2,x_3,\ldots,x_K$ denote the other variables in this constraint. Recall that  the labelling
of variables in this constraint was done consistently with the labelling of the tree
$T$, and, as a result, it
is independent from the measure ${\cal P}(d)$. Let $\sigma:\{1,2,\ldots,K\}\rightarrow\{1,2,\ldots,K\}$
specify the matching between the order in the constraint and in the tree. That is, our
constraint $C(x_1)$ is $\sum_{1\leq k\leq K}a_{rk}x_{\sigma(k)}\leq b_r+W_j+\psi_j$. From
(\ref{eq:OptPsi}) we obtain that, for the limiting probability
measure ${\cal P}(d)$, we also have $\Psi_j=\max(0,\sum_{k}a_{rk}X_{\sigma(k)}-b_r-W_j)$ almost surely.
Let $\bar X_k\equiv \mathbb{E}[X_k|T^W(x_k,d)]$ and let
\be{eq:barPsi}
\bar\Psi_j\equiv\max(0,\sum_{1\leq k\leq K}a_{rk}\bar X_{\sigma(k)}-b_r-W_j).
\ee
That is $\bar\Psi_j$ is the optimal value to be assigned to $\psi_j$ when the variables $x_k$ take values
$\bar X_k$. Observe that, just like $\bar X_k$, $\bar\Psi_j$ is a random variable. Its
value is determined by the random tree $T(x_k,d)$ and the values of $W_j$, and, importantly,
it is different from $\Psi_j$. The following proposition  is the main technical result of this section.
Jumping ahead, the usefulness of this proposition is that if we assign the value $\bar X_k$ to each variable $x_k$ and
take $\psi_j=\bar\Psi_j$ for each constraint containing $x_1$, then we obtain that the corresponding
objective value of the linear program (\ref{eq:GLP}) "per variable" $x_1$ almost
achieves value $\lambda(c)$, provided $d$ is sufficiently large. We will use this in the following
section for the projection step, where for a random instance of linear program (\ref{eq:GLP}),
we assign value $\bar X_i$ for every variable $x_i, 1\leq i\leq n$ and obtain an objective
value close to $\lambda(c)n$, for arbitrary large $n$.

\begin{prop}\label{prop:consistency}
For every $\epsilon>0$ there exists sufficiently large $d=d(\epsilon)$ such that
$\mathbb{E}[w_x\bar X_1+{w_{\psi}\over K}\sum_j\bar\Psi_j]<\lambda(c)+\epsilon$,
where the summation is over all the constraints in ${\cal J}(x_1)$.
\end{prop}

\begin{proof} Fix $\delta>0$ very small. Fix any constraint $C_j\in{\cal J}(x_1)$  and a variable
$x_k$ in it. We first show that for sufficiently large $d$, $\mathbb{P}\{|\bar X_k-\mathbb{E}[X_k|T^W(x_1,d)]|>\delta\}<\delta$.
In other words, the expected values of $X_k$ conditioned on  $T^W(x_k,d)$ and $T^W(x_1,d)$ are sufficiently close to
each other. For this purpose we use the martingale convergence theorem and Proposition \ref{prop:martingale1}.
Take largest integer $t$ such that $tK+k\leq d$. Since, by the first part of the proposition,
$\mathbb{E}[X_k|T^W(x_k,d)], d=1,2,\ldots$ is a martingale,  and since $d-(tK+k)<K$ then, applying
Theorem \ref{theorem:MartingaleConvergence}, $|\mathbb{E}[X_k|T^W(x_k,d)]-\mathbb{E}[X_k|T^W(x_k,tK+k)]|$
becomes very small w.h.p. as $d$ becomes large. In particular,
\[
\mathbb{P}\{|\mathbb{E}[X_k|T^W(x_k,d)]-\mathbb{E}[X_k|T^W(x_k,tK+k)]|>\delta\}<\delta,
\]
for sufficiently large $d$. By the second part of Proposition  \ref{prop:martingale1},
$\mathbb{E}[X_k|T^W(x_{d\mod(K)},d)], d=1,2,\ldots$ is a also a martingale, and again
applying Theorem \ref{theorem:MartingaleConvergence}, we obtain that
\[
\mathbb{P}\{|\mathbb{E}[X_k|T^W(x_k,tK+k)]-\mathbb{E}[X_k|T^W(x_1,tK+1)]|>\delta\}<\delta,
\]
Finally, again applying the first part of Proposition \ref{prop:martingale1} to the
sequence $\mathbb{E}[X_1|T^W(x_1,d)],d=1,2,\ldots$, and using $d-(tK+1)<2K$, we obtain
\[
\mathbb{P}\{|\mathbb{E}[X_k|T^W(x_1,d)]-\mathbb{E}[X_k|T^W(x_1,tK+1)]|>\delta\}<\delta,
\]
for sufficiently large $d$. Combining and using $\bar X_k=\mathbb{E}[X_k|T^W(x_k,d)]$, we obtain
\be{eq:barXk}
\mathbb{P}\{|\bar X_k-\mathbb{E}[X_k|T^W(x_1,d)]|>3\delta\}<\delta,
\ee
as claimed.

For every  constraint $C_j\in{\cal J}(x_1)$ introduce
\be{eq:hatPsi}
\hat \Psi_j\equiv \max(0,\sum_{k}a_{rk}\mathbb{E}[X_{\sigma(k)}|T^W(x_1,d)]-b_r-W_j),
\ee
where $\sum_ka_{rk}x_{\sigma(k)}\leq b_r+W_j+\psi_j$ is the expanded form of the constraint $C_j$.
That is $\hat\Psi_j$ is defined just like $\bar\Psi_j$, except for conditioning
is done on the tree $T^W(x_1,d)$ instead of $T^W(x_k,d)$.
Applying (\ref{eq:barPsi}), (\ref{eq:barXk})
and recalling $\hat\Psi_j,\bar\Psi_j\in[0,B_{\psi}]$, we obtain that for our constraint $C_j$
\be{eq:boundHatPsi}
\mathbb{P}\{|\hat \Psi_j-\bar \Psi_j|>3\delta K\max_{rk}|a_{rk}|\}<K\delta.
\ee
Recall, that, according the measure ${\cal P}(d)$, $J(x_1)$ is distributed as $\Pois(cK)$. Select
a value $M(\delta)>0$ sufficiently large, so that $\mathbb{P}\{J(x_1)>M(\delta)\}<\delta$.
Conditioning on the event $J(x_1)\leq M(\delta)$ we obtain
\begin{align}
&\mathbb{P}\Big\{\sum_j|\hat \Psi_j-\bar \Psi_j|>3\delta M(\delta)K\max_{rk}|a_{rk}|\Big|J(x_1)\leq M(\delta)\Big\} & \nonumber \\
&\leq \mathbb{P}\Big\{\exists C_j\in{\cal J}(x_1):|\hat \Psi_j-\bar \Psi_j|>3\delta K\max_{rk}|a_{rk}|\Big|J(x_1)\leq M(\delta)\Big\} & \nonumber \\
&\leq  M(\delta)\mathbb{P}\Big\{|\hat \Psi_j-\bar \Psi_j|>3\delta K\max_{rk}|a_{rk}|\Big|J(x_1)\leq M(\delta)\Big\} & \label{eq:somej} \\
&\leq  {M(\delta)\mathbb{P}\Big\{|\hat \Psi_j-\bar \Psi_j|>3\delta K\max_{rk}|a_{rk}|\Big\}\over\mathbb{P}\{J(x_1)\leq M(\delta)\}} & \nonumber \\
&\leq  {M(\delta)K\delta\over 1-\delta}<M(\delta)(K+1)\delta, & \nonumber
\end{align}
where the summation is over constraints in ${\cal J}(x_1)$,  the last inequality holds provided $\delta<1/K$,
and $j$ in (\ref{eq:somej})  corresponds to any fixed constraint in ${\cal J}(x_1)$.
Combining with the event $J(x_1)>M$, we obtain from above that, without any conditioning,
\[
\mathbb{P}\{\sum_j|\hat \Psi_j-\bar \Psi_j|>3\delta M(\delta)K\max_{rk}|a_{rk}|\}<M(\delta)(K+1)\delta+\delta<M(\delta)(K+2)\delta.
\]
Since $\hat\Psi_j,\bar\Psi_j\in [0,B_{\psi}]$, then the bound above implies
\be{eq:expHatPsi}
\Big|\mathbb{E}[\sum_j\hat\Psi_j]-\mathbb{E}[\sum_j\bar\Psi_j]\Big|\leq
B_{\psi}M(\delta)(K+2)\delta+3\delta M(\delta)K\max_{rk}|a_{rk}|.
\ee
Our final step is to relate $\hat\Psi_j$ to the random variables $\Psi_j$, where $(X(d),\Psi(d))$ are drawn according
to the probability distribution ${\cal P}(d)$. The distinction between $\hat\Psi_j$ and $\Psi_j$ is somewhat subtle and we expand on
it here. As we observed above, the values of $\Psi(d)$ are determined almost surely by $\Psi_j=\max(0,\sum_ka_{rk}X_{\sigma(k)}-b_r-W_j)$,
with respect to the measure ${\cal P}(d)$, when the corresponding
constraint is $\sum_ka_{rk}x_{\sigma(k)}\leq b_r+W_j+\psi_j$. These values of $\Psi_j$ are different, however, from
$\hat \Psi_j$ which are obtained by first taking the expectations of $X_k$ conditioned on trees $T^W(x_1,d)$ and
then setting $\hat\Psi_j$ according to (\ref{eq:hatPsi}). Naturally, $\Psi_j$ and $\hat\Psi_j$ are related to
each other. Note, that for every constraint $C_j$ almost surely
\be{eq:lowerIneq}
\Psi_j\geq\sum_ka_{rk}X_{\sigma(k)}-b_r-W_j, \qquad \Psi_j\geq 0.
\ee
By taking the conditional expectations and using the linearity of inequalities above, we obtain
\be{eq:lowerIneqExp}
\mathbb{E}[\Psi_j|T^W(x_1,d)]\geq\sum_ka_{rk}\mathbb{E}[X_{\sigma(k)}|T^W(x_1,d)]-b_r-W_j, \qquad \mathbb{E}[\Psi_j|T^W(x_1,d)]\geq 0,
\ee
where we use the trivial equality $\mathbb{E}[W_j|T^W(x_1,d)]=W_j$.
From the definition of $\hat\Psi_j$ in (\ref{eq:hatPsi}), we obtain then
$\mathbb{E}[\Psi_j|T^W(x_1,d)]\geq\hat\Psi_j$ almost surely, with respect to the random variables $T^W(x_1,d)$.
As a result $\sum_j\mathbb{E}[\Psi_j|T^W(x_1,d)]\geq \sum_j\hat\Psi_j$, where the summation is over constraints in ${\cal J}(x_1)$.
Therefore,
\be{eq:E>hatE}
\mathbb{E}[\sum_j\Psi_j]\geq\mathbb{E}[\sum_j\hat\Psi_j].
\ee
Recall from the last part of Proposition \ref{prop:ConvergenceMeasures}, that with respect to measure ${\cal P}(d)$,
\[
\mathbb{E}[w_xX_1+{w_{\psi}\over K}\sum_j\Psi_j]=\lambda(c).
\]
Combining this with (\ref{eq:expHatPsi}), (\ref{eq:E>hatE}), and using a simple observation
$\mathbb{E}[\bar X_1]=\mathbb{E}[\mathbb{E}[X_1|T^W(x_1,d)]]=\mathbb{E}[X_1]$, we obtain
\begin{align}
&\mathbb{E}[w_x\bar X_1+{w_{\psi}\over K}\sum_j\bar\Psi_j]  \nonumber \\
&\leq \mathbb{E}[w_xX_1+{w_{\psi}\over K}\sum_j\Psi_j]+w_{\psi}B_{\psi}M(\delta)(K+2)\delta+3\delta M(\delta)K\max_{rk}|a_{rk}| \nonumber \\
&=\lambda(c)+w_{\psi}B_{\psi}M(\delta)(K+2)\delta+3\delta M(\delta)K\max_{rk}|a_{rk}|. \label{eq:boundFinal}
\end{align}
Recall that $J(x_1)$ is distributed as $\Pois(cK)$ and, in particular, has exponentially decaying tails.
Therefore, for any $\epsilon>0$ we can find sufficiently small $\delta$ and the corresponding
$M(\delta)$ such that $\delta M(\delta)<\epsilon$ and $\mathbb{P}\{J(x_1)>M(\delta)\}<\delta$. All the other
values in the right-hand side of the bound (\ref{eq:boundFinal}) are constants. Therefore, for any $\epsilon>0$,
we can find sufficiently small $\delta>0$ such that the right-hand side is at most $\lambda(c)+\epsilon$. Choosing
$d$ sufficiently large for this $\delta$, we obtain the result. \qed
\end{proof}

\section{Projection}\label{section:projection}
In this section we complete the proof of Theorem \ref{theorem:mainMaxGen} by
proving the existence of the  limit (\ref{eq:MaxLSATEXP}). We use the limiting distribution
${\cal P}(d)$ constructed in Section \ref{sec:PoissonTree} and "project" it onto
any random instance of linear program (\ref{eq:GLP}) with $n$ $x$-variables and $cn$  constraints
$C_1,\ldots,C_{cn}$.

Fix $c>0$ and $\epsilon>0$ and take  $n$ to be a large integer. We construct a feasible solution
$x_i\in [B^1_x,B^2_x], \psi_j\geq 0, 1\leq i\leq n, 1\leq j\leq cn$ as follows.
For each variable $x_i$ consider its depth-$d$ neighborhood ${\cal B}(x_i,d,n)$,
where $d=d(\epsilon)$ is taken as  in Proposition \ref{prop:consistency}.
If the ${\cal B}(x_i,t,n)$ is a tree $T(x_i,d)$, then we set the value of $x_i$ equal to
$X_i(n)\equiv\mathbb{E}[X_i|T^W(x_i,d)]$, where the expectation is with respect to the measure ${\cal P}(d)$.
That is we observe the depth-$d$ tree $T(x_i,d)$ together with values $W_j$ corresponding
to the constraints $C_j$ in this tree and set the values $x_i$ to be equal to the expectation of $X_i$ condition on this observation.
If, on the other hand, ${\cal B}(x_i,t,n)$ is not a tree, then we assign any value to $x_i=X_i(n)$,
for example $x_i=B^1_x$. Once we have assigned values to $x_i$ in the manner above, for every constraint
$C_j:\sum_k a_{rk}x_{i_k}\leq b_r+W_j+\psi_j$, we set its corresponding value of $\psi_j$ to
$\Psi_j(n)\equiv \max\{0,\sum_ka_{rk}x_{i_k}-b_r-W_j\}$ - the optimal choice for  given values of $x_i$-s.

\begin{prop}\label{prop:projection}
For every $\epsilon>0$, the solution constructed above has expected cost at most $(\lambda(c)+\epsilon)n$
for all sufficiently large $n$.
\end{prop}

Since $\epsilon$ was an arbitrary constant and since $\lambda(c)$ was defined by (\ref{eq:lambdalimsup}) the proposition shows that the
assignment above satisfies $\lim_{n\rightarrow\infty}$ $\mathbb{E}[{\cal GLP}(n,c)]/n=\lambda(c)$.
Therefore the proposition implies Theorem \ref{theorem:mainMaxGen}.

\begin{proof}
Select one of the $n$ variables $x_i$ uniformly and random. W.l.g. assume it is $x_1$. Fix $\epsilon_0>0$ very small.
We claim that when $n$ is sufficiently large, we have
\be{eq:boundEpsilon0}
\mathbb{E}[w_xX_1(n)+{w_{\psi}\over K}\sum_j\Psi_j(n)]\leq \lambda(c)+3\epsilon_0,
\ee
where the summation is over all constraints $C_j$ containing $x_1$.
First let us show that this implies the statement of the proposition.
Indeed, multiplying the inequality by $n$, and recalling that $x_1$ was uniformly selected from $x_i,1\leq i\leq n$, we
obtain that $\mathbb{E}[w_x\sum_{1\leq i\leq n}X_i(n)+w_{\psi}\sum_j\Psi_j(n)]\leq (\lambda(c)+3\epsilon_0)n$,
where again we used the fact that each variable $\Psi_j(n)$ was counted exactly $K$ times. By taking
$\epsilon_0<\epsilon/3$,  we obtain the required bound.

Consider the neighborhood ${\cal B}(x_1,d+1,n)$ and suppose first that it is not a tree. Denote this event
by ${\cal E}_n$. From the second part of Lemma \ref{lemma:cycles}, $\mathbb{P}\{{\cal E}_n\}=O(1/n)$,
where the notation $O(\cdot)$ involves constants $d,c$ and $K$. Recall, that the values of $\Psi_j(n)$
never exceed $B_{\psi}$. Then we have
\begin{align}
&\mathbb{E}[w_xX_1(n)+{w_{\psi}\over K}\sum_j\Psi_j(n)|{\cal E}_n] \nonumber \\
&\leq w_x\max\{|B^1_x|,|B^2_x|\}+w_{\psi}B_{\psi}\mathbb{E}[|J(x_1,n)|{\cal E}_n] \nonumber \\
&=w_x\max\{|B^1_x|,|B^2_x|\}+w_{\psi}B_{\psi}{\mathbb{E}[|J(x_1,n)1\{{\cal E}_n\}]\over \mathbb{P}\{{\cal E}_n\}}, \nonumber
\end{align}
where $1\{{\cal E}_n\}$ is the indicator of an event ${\cal E}_n$. Recall from the proof of Lemma \ref{lemma:DegreePoisson},
that the probability $J(x_1,n)=s$ is asymptotically $\Big(\begin{array}{c} cn \\s\end{array}\Big)
(K/n)^s(1-K/n)^{cn-s}$ for large $n$. It follows that the probability that $J(x_1,n)$ exceeds $\log^2n$ is at most
$1/n^3$ for sufficiently large $n$ (more accurate bounds can be obtained \cite{JansonBook}, which are not required here).
Then
\begin{align}
&\mathbb{E}[J(x_1,n)1\{{\cal E}_n\}] \nonumber \\
&\leq\mathbb{E}[J(x_1,n)1\{{\cal E}_n\}1\{J(x_1,n)\leq \log^2n\}]+\mathbb{E}[J(x_1,n)1\{{\cal E}_n\}1\{J(x_1,n)>\log^2n\}] \nonumber \\
&\leq\log^2n\mathbb{P}\{{\cal E}_n\}+{cn\over n^3}\mathbb{P}\{{\cal E}_n\},
\end{align}
where we used the fact that $J(x_1,n)\leq cn$ with probability one. But since $\mathbb{P}\{{\cal E}_n\}=O(1/n)$, we obtain
that $\mathbb{E}[(w_xX_1(n)+{w_{\psi}\over K}\sum_j\Psi_j(n)) 1\{{\cal E}_n\}]=o(1)$.

Suppose now that ${\cal B}(x_1,d+1,n)$ is a tree $T(x_1,d)$, that is the event $\bar{\cal E}_n$ occurs.
Select any of the constraints containing $x_1$ (if any exist), and let $x_2,\ldots,x_K$ be the variables in this constraints.
Note that then ${\cal B}(x_k,d,n)$ are also trees $T(x_k,d)$ and the corresponding values
$X_k(n)=\mathbb{E}[X_k|T^W(x_k,d)], k=2,3,\ldots,K$ for these variables
are set based only on the observed trees $T(x_k,d)$ and values $W_j$, that is exactly
as $X_k$ where defined in the previous section. Then the same correspondence holds between
$\Psi_j(n)$ and $\bar\Psi_j$. Applying Proposition \ref{prop:consistency}, and the first
part of Proposition \ref{prop:ConvergenceMeasures}, that is the fact that ${\cal B}(x_1,n,d)$
converges to $T(x_1,d)$ distributed according to ${\cal P}(d)$, we obtain that, for sufficiently
large $n$, $\mathbb{E}[w_xX_1(n)+{w_{\psi}\over K}\sum_j\Psi_j(n)]\leq \lambda(c)+2\epsilon_0$,
where the second $\epsilon_0$ comes from approximating ${\cal B}(x_1,n,d)$ by $T(x_1,d)$.
Combining with the case of non-tree ${\cal B}(x_1,d,n)$ we obtain
\[
\mathbb{E}[w_xX_1(n)+{w_{\psi}\over K}\sum_j\Psi_j(n)]\leq \lambda(c)+2\epsilon_0+o(1)<\lambda(c)+3\epsilon_0.
\]
for sufficiently large $n$, just as required by (\ref{eq:boundEpsilon0}). This concludes the proof of Theorem \ref{theorem:mainMaxGen}. \qed
\end{proof}

\section{Applications to  maximum weight $b$-matchings in sparse random graphs}\label{section:MatchingCyclePacking}
The main goal of this section is proving Theorem \ref{theorem:mainMatching}.
We begin with  a linear programming formulation of the maximum weight matching problem. Suppose we have (a non-random) graph
with $n$ nodes and $m$ undirected edges represented as pairs $(i,j)$ of nodes. Denote by $E$ the edge set of the graph.
The edges are equipped with (non-random) weights
$0\leq w_{i,j}\leq w_{\rm max}$. Given  $V\subset [n]$, let $\delta(V)$ denote the set of edges with exactly one
end point in $V$. A classical result from the theory of combinatorial optimization (\cite{SchrijverBook}, Theorem 32.2) states that the following
linear programming problem provides an exact solution (namely, it  is a tight relaxation) of the maximum weight
$b$-matching problem:
\begin{align}
&{\rm Maximize}\,\,\sum_{i,j}w_{i,j}x_{i,j}   \label{eq:MatchingObjective}\\
&{\rm subject\,\,to:} & \nonumber \\
&\sum_{j}x_{i,j}\leq b, \qquad \forall i=1,2,\ldots,n,  \label{eq:MatchingDegree} \\
&\sum_{i,j\in V}x_{i,j}+\sum_{(i,j)\in A}x_{i,j}\leq {b|V|+|A|-1\over 2},
\qquad \forall ~V\subset[n], A\subset \delta(V) {\rm ~such~that}~b|V|+|A|{\rm~is~odd},  \label{eq:MatchingOddCycle} \\
&0\leq x_{i,j}\leq 1 \label{eq:nonnegative}.
\end{align}

Specifically, there exists an optimal solution of this linear programming problem which is always integral and it
corresponds to the maximum weight $b$-matching. We denote the optimal value of the linear program above by ${\cal LPM}(G)$.

Our plan for proving Theorem \ref{theorem:mainMatching} is as follows. We first show that when
the graph has very few small cycles  (and this will turn out to be the case for $G(n, cn))$,
the optimal value ${\cal LPM}^0(G)$ of the modified linear program, obtained from (\ref{eq:MatchingObjective})-(\ref{eq:nonnegative})
by dropping the constraints (\ref{eq:MatchingOddCycle}), is very close to ${\cal LPM}(G)$. In the context of random graphs
$G(n, cn)$, this will imply that the difference
$\Big|{\cal LPM}(G(n, cn))-{\cal LPM}^0(G(n, cn))\Big|=o(n)$, w.h.p. We then take
the dual of the modified linear program (\ref{eq:MatchingObjective}), (\ref{eq:MatchingDegree}), (\ref{eq:nonnegative})
and show that it has the form (\ref{eq:GLP}). Applying Theorem \ref{theorem:mainMaxGen} we will obtain that the limit
$\lim_n{\cal LPM}^0(G(n, cn)/n$ exists w.h.p. This will imply the existence of the limit
$\lim_n{\cal LPM}(G(n, cn)/n$ w.h.p. Since ${\cal LPM}(G(n, cn)$ is the maximum
weight $b$-matching in $G(n, cn)$, that is ${\cal M}_w(n,c,b)$, and Theorem \ref{theorem:mainMatching}  follows.

Naturally, ${\cal LPM}^0(G)\geq {\cal LPM}(G)$.

\begin{prop}\label{prop:RelaxedLP}
Given a weighted graph $G$ and $d\geq 3$ let $L(d)$ denote the set of cycles of length
$<d$ in $G$ and let $M(d)$ denote the total number of edges in $L(d)$. Then
\be{eq:RelaxedLP}
{d-1\over d}({\cal LPM}^0(G)-M(d)w_{\max})\leq {\cal LPM}(G)\leq {\cal LPM}^0(G).
\ee
\end{prop}

\begin{proof} We already observed that the right-hand side of (\ref{eq:RelaxedLP}) holds. We now concentrate
on the left-hand side bound.

Let $x^0=(x^0)_{i,j}$ be an optimal solution of the  linear program (\ref{eq:MatchingObjective}),
(\ref{eq:MatchingDegree}), (\ref{eq:nonnegative}) with its the optimal value
${\cal LPM}^0(G)$. In the graph $G$ we delete all the $M(d)$ edges which contribute to $L(d)$ together
with the corresponding values of $x^0_{i,j}$. Consider the resulting solution $x^1=(x^1_{i,j})$ to the linear program
corresponding to the reduced weighted graph $G^1$, which now does not contain any cycles of length less
than $d$. The objective value of the linear program (\ref{eq:MatchingObjective}),
(\ref{eq:MatchingDegree}), (\ref{eq:nonnegative}) corresponding to the solution $x^1$ is at least
${\cal LPM}^0(G)-M(d)w_{\max}$, since,  by constraint (\ref{eq:nonnegative}), $x^0_{i,j}\leq 1$
for all the edges $(i,j)$. We further modify the solution $x^1$ to $x^2$ by letting $x^2_{i,j}=(1-1/d)x^1_{i,j}$
for every edge $(i,j)$ in the graph $G^1$. The objective value of the linear program (\ref{eq:MatchingObjective}),
(\ref{eq:MatchingDegree}), (\ref{eq:nonnegative}) corresponding to $x^2$ is then at least
${d-1\over d}({\cal LPM}^0(G)-M(d)w_{\max})$. We claim that in fact $x^2$ is a feasible solution
to the linear program (\ref{eq:MatchingObjective})-(\ref{eq:nonnegative}), implying
(\ref{eq:RelaxedLP}) and completing the proof of the proposition.

Clearly, constraints (\ref{eq:MatchingDegree}) and (\ref{eq:nonnegative}) still hold. We concentrate on
(\ref{eq:MatchingDegree}). Consider any  set $V\subset[n]$ and $A\subset \delta(V)$ such that $b|V|+|A|$ is odd.
Assume first $|V|+|A|<d$. Let $\hat V$ denote the union of $V$ and the end points of edges in $A$.
Then $|\hat V|<d$. Since $G^1$ does not contain any cycles of size $<d$, then $\hat V$ does not contain any cycles at all,
and therefore is a forrest. In particular it is a bipartite graph. Let $\hat x^1$ denote the sub-vector of the
vector $x^1$ corresponding to edges with both ends $\hat V$. Since  (\ref{eq:MatchingDegree}) holds for $x^1$, then it also
holds for the vector $\hat x^1$ for all nodes in $\hat V$. A classical result from a combinatorial
optimization theory \cite{SchrijverBook} states that for every bipartite graph, the polytop corresponding to the
degree constraints (\ref{eq:MatchingDegree})  has only integral extreme points, which are $b$-matchings
(in the reduced graph $\hat V$). This follows from the fact that this polytop when described in matrix form
$Bx\leq b$ corresponds to the case when $B$ is totally unimodular. We refer the reader to \cite{SchrijverBook} for
the details. But since any integral solution $\hat x^1$ corresponding to the $b$-matching must satisfy the
constraints (\ref{eq:MatchingOddCycle}) by Theorem 32.2 in \cite{SchrijverBook},
then these constraints are automatically satisfied by $x^1$. Moreover then their are satisfied by $x^2$.
We proved that (\ref{eq:MatchingOddCycle}) holds whenever $|V|+|A|<d$.

Suppose now $|V|+|A|\geq d$. For the solution $x^1$, let us sum the constraints (\ref{eq:MatchingDegree}) corresponding to all
the nodes in $V$ and sum the right-hand side constraints (\ref{eq:nonnegative}) corresponding to all the
edges in $A$. Each value $x_{i,j}$ for $i,j\in V$ is counted twice, once for node $i$ and once for node $j$.
Each value $x_{i,j}$ for $(i,j)\in A$ is also counted twice, once for constraint (\ref{eq:MatchingDegree}) for
the node $i$ and once for constraint (\ref{eq:nonnegative}) for the edge $(i,j)$. Then we obtain
$2(\sum_{i,j\in V}x^1_{i,j}+\sum_{(i,j)\in A}x^1_{i,j})\leq b|V|+|A|$. This implies
$\sum_{i,j\in V}x^2_{i,j}+\sum_{(i,j)\in A}x^2_{i,j}\leq (1-1/d){b|V|+|A|\over 2}\leq {b|V|+|A|-1\over 2}$,
since by assumption, $b|V|+|A|\geq |V|+|A|\geq d$. Again we showed that the constraint (\ref{eq:MatchingDegree}) holds for the
solution $x^2$. \qed
\end{proof}

We return to our main setting -- sparse random graph $G(n, cn)$, with edges equipped with randomly
generated weights $W_{i,j}$, drawn according to a distribution function $w^{\rm edge}(t), t\geq 0$ with support
$[0,B_w]$. Let $E=E(G(n, cn))$ denote the edge set of this graph. We denote the value of  ${\cal LPM}^0(G(n, cn))$ by
${\cal LPM}^0(n,c)$ for simplicity. That is, ${\cal LPM}^0(n,c)$ the optimal (random) value of the linear program
(\ref{eq:MatchingObjective}), (\ref{eq:MatchingDegree}), (\ref{eq:nonnegative}) on the graph $G(n, cn)$.

\begin{prop}\label{prop:MatchingGapSparseGraphs}
W.h.p. as $n\rightarrow\infty$
\be{eq:RelaxedLPRandom}
{\cal LPM}^0(n,c)-o(n)\leq {\cal M}_w(n,c,b)\leq {\cal LPM}^0(n,c).
\ee
\end{prop}

\begin{proof}
Let $d(n)$ be a very slowly growing function of $n$. It is well known that in $G(n, cn)$
w.h.p. the total number of edges which belong to at least one cycle with size $<d(n)$ is $o(n)$ (far more
accurate bounds can be obtained \cite{JansonBook}). Thus $M(d(n))=o(n)$ w.h.p. Note also that from
(\ref{eq:nonnegative}), we have ${\cal LPM}^0(n,c)\leq w_{\max}bn=O(n)$. Applying Proposition
\ref{prop:RelaxedLP} with $d=d(n)$, we obtain the result. \qed
\end{proof}

Our final goal is proving the convergence w.h.p. of ${\cal LPM}^0(n,c)/n$. We use linear programming duality for this purpose.
Consider the  dual of the  linear program (\ref{eq:MatchingObjective}), (\ref{eq:MatchingDegree}), (\ref{eq:nonnegative}),
generated on the weighted graph $G(n, cn)$. It involves variables $y_1,\ldots,y_n$ and has the following form.

\begin{align}
&{\rm Minimize}\,\,b\sum_{1\leq i\leq n}y_i+\sum_{(i,j)\in E}\psi_{i,j}   \label{eq:MatchingDualObjective}\\
&{\rm subject\,\,to:} & \nonumber \\
&y_i+y_j+\psi_{i,j}\geq W_{i,j}, \qquad \forall~(i,j)\in E  \label{eq:MatchingDualConstraint} \\
&y_{i},\psi_{i,j}\geq 0 \label{eq:DualNonnegative}.
\end{align}

The objective value of this linear program is also ${\cal LPM}^0(n,c)$, thanks to  the strong duality of linear programming.
The linear program above is almost of the form (\ref{eq:GLP}) that we need in order to apply Theorem \ref{theorem:mainMaxGen}.
Let us rewrite the linear programm above in the following equivalent form
\begin{align}
&{\rm Minimize}\,\,b\sum_{1\leq i\leq n}y_i+\sum_{(i,j)\in E}\psi_{i,j}   \label{eq:MatchingDualObjective1}\\
&{\rm subject\,\,to:} & \nonumber \\
&(-1)y_i+(-1)y_j\leq -W_{i,j}+\psi_{i,j}, \qquad \forall~(i,j)\in E  \label{eq:MatchingDualConstraint1} \\
&y_{i},\psi_{i,j}\geq 0 \label{eq:DualNonnegative1},
\end{align}

Let the set of constraints ${\cal C}$ of the linear program (\ref{eq:GLP}) contain only one
element $(-1)y_1+(-1)y_2\leq 0$. We set $B^1_x=0,B^2_x=B_w$, where, as we recall, $[0,B_w]$ is the support of the distribution of $W_{i,j}$.
We also set $w_x=b$ and $w_{\psi}=1$.  Our linear program has now the form (\ref{eq:GLP}) except for we
need to consider in addition the constraints $y_i\leq B_w$. We claim that in fact
these constraints are redundant. In fact any value of $y_i$ which exceeds $B_w\geq W_{i,j}$ can be decreased to $B_w$,
resulting in a smaller value of the  objective function and still honoring all the constraints. Thus, we may replace
$0\leq y_i\leq B_w$ simply by $y_i\geq 0$. We conclude that the linear program
(\ref{eq:MatchingDualObjective1}), (\ref{eq:MatchingDualConstraint1}), (\ref{eq:DualNonnegative1})
has form (\ref{eq:GLP}). Applying Theorem \ref{theorem:mainMaxGen}, there exists a function $g(c)\geq 0$ such that w.h.p.
${\cal LPM}^0(n,c)/n\rightarrow g(c)$ as $n\rightarrow\infty$. Finally, applying Proposition \ref{prop:MatchingGapSparseGraphs},
we obtain (\ref{eq:LimitMatching}). This concludes the proof of Theorem \ref{theorem:mainMatching}. \qed

\section{Discussion}\label{section:discussion}
The results of the present paper lead to  several interesting open questions.
In addition to Conjecture \ref{conj:LSAT}, stated in Section \ref{section:model}, it seems that the following analogue of
Conjecture \ref{conj:MaxSAT} is reasonable.

\begin{conj}\label{conj:MaxLSAT}
Consider random K-LSAT problem (linear program (\ref{eq:LP})) with $n$ variables and $m$ constraints, where all
$\psi$ variables are set to zero. Let $N(n,m)$ denote the maximum cardinality subset of constraints $C_j$ which
is feasible. For any $c>0$, the limit $\lim_{n\rightarrow\infty}(N(n,cn)/n)$ exists.
Moreover, this limit is strictly smaller than unity, for all $c>c^*_K$.
\end{conj}


An interesting group of questions relates to the behavior of the function
$f(c)=\lim_n \mathbb{E}[{\cal GLP}(n)]/n$, which, by results of this paper is equal to zero for $c<c^*$,
in the specific context of linear program (\ref{eq:LP}).
What can be said about $f(c)$ near $c^*_K$, when $c^*_K<\infty$?
It is not hard to show that $f(c)$ is continuous and non-decreasing function of $c$. Is it differentiable
in $c^*_K$? Is it convex, concave or neither? Similar questions arise in connection with percolation probability $\theta(p)$
which is a function of a bond or site probability $p$. It is known that for every
dimension $d$ there exists a critical value $p^*_d$, for which
the probability $\theta(p)$ of existence of an infinite component containing
the origin in $d$-dimensional bond/site percolation model, is equal to zero for $p<p^*_d$ and
is positive for $p>p^*_d$, (whether $\theta(p)=0$ for $p=p^*_d$ is a major outstanding problem
in percolation theory), see \cite{grimmett}.
The behavior of $\theta(p)$ near $p^*_d$ is a well-known open problem. We also refer the
reader to \cite{AldousPercus} for the related discussion on scaling limits and universality
of random combinatorial problems.

Computing the values of $c^*_K$ is another question which seems so far beyond the techniques
applied in this paper. The cases where local weak convergence methods lead to computation of
limiting values seem to be related  to special structures of the corresponding problems.
For example a very clever recursion was used by Aldous to prove $\zeta(2)=\pi^2/6$ limit for the
expected value of the minimum weight assignment, \cite{Aldous:assignment00}. See also Bandyopadhyay \cite{Bandyopadhyay}, who
investigates some properties of the limiting measure arising in the proof of $\zeta(2)$ result.

Similarly, a special
structure was used in \cite{AldousSteele:survey} to prove $\zeta(3)$ limit for
minimum spanning tree on a complete graph with exponentially distributed weights.
(This problem was originally solved by Frieze \cite{frieze}, using combinatorial methods).
Uncovering such special structure for our problem for the purposes of computing $c^*_K$
is an interesting question.

A separate course of investigation is to formulate and study randomly generated integer programming
problems as a unifying framework for studying random k-SAT, coloring, maximum k-cuts , maximum independent sets,
and other problems. For example, under which conditions on the set of prototype constraints does the feasibility
problem experience a sharp transition? These conditions should be generic enough to include the problems
mentioned above. Also we suspect that other results from polydedral combinatorics can be of use when
studying these problems within an integer programming framework.

Finally, it should be clear  what
goes wrong when one tries to use local weak convergence approach for random K-SAT problem, for example
along the lines of Theorem \ref{theorem:mainMaxGen}. Our approach is built on using the values like
 $\mathbb{E}[X_k|\,\cdot\,]$ to construct a feasible solution, but these expectations are not necessarily integers.
Digging somewhat deeper into the issue, it seems that local
weak convergence method in general is not very hopeful for resolving Conjecture \ref{conj:main},
since it looks only into constant size neighborhoods of nodes. To elaborate somewhat this point
consider maximal independent set problem in $r$-regular random graphs, discussed in \cite{AldousSteele:survey}.
For almost any node in such a graph, its constant size neighborhood is a $r$-regular tree, and,
as such, the neighborhoods are indistinguishable.
In such circumstances it seems hard to try to
concoct a solution which is based only on neighborhoods of nodes. Some long-range structural properties
of these graphs like structure of cycles have to be considered. We refer the reader to
\cite{AldousSteele:survey} for a further discussion of this issue.

\section*{Acknowledgments} The author gratefully acknowledges interesting and fruitful discussions with
David Aldous, Michael Steele and Alan Frieze. The author is grateful to Jonathan Lee for references to
the results in polyhedral combinatorics, and  to Maxim Sviridenko for suggesting the proof of Proposition \ref{prop:RelaxedLP}.

\bibliographystyle{amsalpha}

\providecommand{\bysame}{\leavevmode\hbox
to3em{\hrulefill}\thinspace}
\providecommand{\MR}{\relax\ifhmode\unskip\space\fi MR }
\providecommand{\MRhref}[2]{%
  \href{http://www.ams.org/mathscinet-getitem?mr=#1}{#2}
} \providecommand{\href}[2]{#2}

\end{document}